\documentclass{article}

\usepackage{graphicx}

\usepackage{latexsym}
\usepackage{amsfonts}
\usepackage[psamsfonts]{amssymb}
\usepackage{enumerate}

\RequirePackage[OT1]{fontenc}
\RequirePackage[amsthm,amsmath]{imsart}

\theoremstyle{plain} 
\newtheorem{theorem}{Theorem}[section]
\newtheorem{corollary}[theorem]{Corollary}
\newtheorem{lemma}{Lemma}[subsection]
\newtheorem{proposition}[theorem]{Proposition}

\theoremstyle{definition} 
\newtheorem{definition}[theorem]{Definition}
\theoremstyle{definition} 

\theoremstyle{remark} 

\theoremstyle{remark} 
\newtheorem{remark}[theorem]{Remark}
\newtheorem*{remark*}{Remark}
\numberwithin{equation}{section}
 
\renewcommand{\r}{\mathsf{r}}
\newcommand{\sign}{\operatorname{sign}}

\newcommand{\esssup}{\operatorname{ess\,sup}}

\newcommand{\lp}{\left(}
\newcommand{\rp}{\right)}

\renewcommand{\le}{\leqslant}
\renewcommand{\ge}{\geqslant}
\renewcommand{\leq}{\leqslant}
\renewcommand{\geq}{\geqslant}
\renewcommand{\succeq}{\succcurlyeq}

\renewcommand{\smallint}{\textstyle{\int}}

\newcommand{\lc}{\mathsf{L\!C}}
\newcommand{\lin}{\mathsf{Lin}}

\renewcommand{\aa}{\overline{a}}
\newcommand{\xx}{\mathbf{x}}

\newcommand{\al}{\alpha}

\renewcommand{\th}{\theta}
\newcommand{\g}{\gamma}
\newcommand{\si}{\sigma}
\newcommand{\la}{\lambda}
\newcommand{\de}{\delta}
\newcommand{\De}{\Delta}

\renewcommand{\Psi}{\overline{\Phi}}

\newcommand{\iid}{\overset{\mathrm{i.i.d.}}{\sim}}

\newcommand{\ii}[1]{\boldsymbol{I}\!\left\{#1\right\}}

\newcommand{\st}{\mathsf{st}}
\newcommand{\ST}{\mathsf{ST}}
\newcommand{\SST}{\mathrm{S\!T}}

\newcommand{\bs}{\mathsf{BS}}
\newcommand{\BS}{\mathrm{B\!S}}
\newcommand{\BC}{\mathrm{B\!C}}
 
\newcommand{\kurt}{\operatorname{\mathrm{kurt}}} 

\renewcommand{\P}{\operatorname{\mathsf{P}}} 
\newcommand{\PP}{\operatorname{\mathsf{P}}}
\newcommand{\E}{\operatorname{\mathsf{E}}}
\newcommand{\Var}{\operatorname{\mathsf{Var}}}

\newcommand{\Z}{\mathbb{Z}}
\newcommand{\R}{\mathbb{R}}
\newcommand{\C}{\mathcal{C}}
\newcommand{\M}[1]{\mathcal{M}_+^{#1}}

\newcommand{\F}[1]{\mathcal{F}_+^{#1}}
\renewcommand{\H}[1]{\mathcal{H}_+^{#1}}

\newcommand{\Fminus}[1]{\mathcal{F}_-^{#1}}
\newcommand{\FF}[1]{\mathcal{F}^{#1}}
\newcommand{\GG}[1]{\mathcal{G}^{#1}}
\newcommand{\G}[1]{\mathcal{G}^{#1}_+}
\newcommand{\GPlus}[1]{\mathcal{G}^{#1}_{++}}

\newcommand{\opt}{\textsf{opt}}

\newcommand{\A}{\mathcal{A}}

\renewcommand{\d}{\mathrm{d}}
\newcommand{\cl}{\operatorname{cl}}

\newcommand{\vp}{\varepsilon}

\begin{document}


\begin{frontmatter}

\title{Exact inequalities for sums\\ of asymmetric random variables,\\ with applications}
\runtitle{Inequalities for sums of asymmetric random variables}

\begin{aug}
\author{\fnms{Iosif} \snm{Pinelis}\ead[label=e1]{ipinelis@math.mtu.edu}}
\runauthor{Iosif Pinelis}


\address{Department of Mathematical Sciences\\
Michigan Technological University\\
Houghton, Michigan 49931, USA\\
E-mail: \printead[ipinelis@mtu.edu]{e1}}
\end{aug}






\begin{abstract}
Let $\BS_1,\dots,\BS_n$ be independent identically distributed random variables each having the standardized Bernoulli distribution with parameter $p\in(0,1)$. Let 
$m_*(p):=(1 + p + 2\,p^2)/(2\sqrt{p - p^2} + 4\,p^2)$ if $0<p\le\frac12$ 
and
$m_*(p):=1$ if $\frac12\le p<1$.  
Let  
$m\ge m_*(p)$.
Let $f$ be such a function that $f$ and $f''$ are nondecreasing and convex.
Then it is proved that for all nonnegative numbers 
$c_1,\dots,c_n$ one has the inequality
\begin{center}
\makebox{
\qquad\qquad
$\E f(c_1\BS_1+\dots+c_n\BS_n)\le\E f\lp s^{(m)}\cdot(\BS_1+\dots+\BS_n)\rp$,
\qquad\quad\ 
}
\end{center}
where 
$s^{(m)}:=\lp\frac1n\,\sum_{i=1}^n c_i^{2m}\rp^\frac1{2m}$.
The lower bound $m_*(p)$ on $m$ is exact for each $p\in(0,1)$.
Moreover, $\E f(c_1\BS_1+\dots+c_n\BS_n)$ is Schur-concave in $(c_1^{2m},\dots,c_n^{2m})$.
A number of related results are presented, including ones for the ``symmetric'' case. 

A number of corollaries are obtained, including upper bounds on generalized moments and tail probabilities of (super)martin\-ga\-les with differences of bounded asymmetry, and also upper bounds on the maximal function of such (super)martin\-ga\-les. It is shown that these results may be important in certain statistical applications. 
\end{abstract}

\begin{keyword}[class=AMS]
\kwd[Primary ]{60E15, 60G50, 60G42, \\
60G48, 62F03, 62F25, 62G10, 60G15}
\kwd[; secondary ]{60E05}
\kwd{62E10}
\kwd{62G35}
\end{keyword}

\begin{keyword}
\kwd{supermartingales}
\kwd{martingales}
\kwd{upper bounds}
\kwd{probability inequalities}
\kwd{generalized moments}
\kwd{$t$ statistic}
\kwd{self-normalized sums}
\end{keyword}





\end{frontmatter}

\tableofcontents

\section{Introduction}\label{intro}
Let $\vp_1,\dots,\vp_n$ be independent Rademacher random variables (r.v.'s), so that $\P(\vp_i=1)=\P(\vp_i=-1)=\frac12$ for all $i$. Let $Z\sim N(0,1)$. 
Let $a_1,\dots,a_n$ be any real numbers such that 
$$a_1^2+\dots+a_n^2=1.$$
The sharp form, 
\begin{equation}\label{eq:khin}
\E f\lp\vp_1a_1+\dots+\vp_n a_n\rp\le\E f(Z),
\end{equation}
of Khinchin's inequality \cite{kh} for $f(x)\equiv |x|^p$ 
was proved by
Whittle (1960) \cite{whittle} for $p\ge3$ and Haagerup (1982) \cite{haag} for $p\ge2$. 

For $f(x)\equiv e^{\la x}$ ($\la\ge0$), inequality \eqref{eq:khin} follows from Hoeffding (1963) \cite{hoeff}, whence
\begin{equation}\label{eq:exp}
\P\lp\vp_1a_1+\dots+\vp_n a_n\ge x\rp\le 
\inf_{\la\ge0}\frac{\E e^{\la Z}}{e^{\la x}}
=e^{-x^2/2}\quad \forall x\ge0.
\end{equation}
Since 
$\P(Z\ge x)\sim\frac1{x\sqrt{2\pi}}e^{-x^2/2}$ as $x\to\infty$,
a factor $\asymp\frac1x$ is ``missing" here. 
The apparent cause of this deficiency is that the class of exponential moment functions $f(x)\equiv e^{\la x}$ ($\la\ge0$) is too small (and so is the class of the power functions $f(x)\equiv |x|^p$). 
 
For all $\al\ge0$, consider the following much richer classes of functions $f\colon\R\to\R$:
\begin{equation}\label{eq:H}
\H\al:=\{f\colon
f(x)=\smallint_{-\infty}^\infty (x-t)_+^\al\,\mu(dt)\quad \forall u\in\R\},
\end{equation}
where $\mu\ge0$ is a Borel measure,
$x_+:=\max(0,x)$ and $x_+^\al:=(x_+)^\al$ for $x\in\R$,
$0^0:=0$; however, the subscript ${}_+$ will have a different meaning when used with functions or classes of functions (as, for example, in the symbol $\H\al$). 

It is easy to see \cite[Proposition 1(ii)]{pin99} that
\begin{equation}\label{eq:F-al-beta}
0\le\beta<\al\quad\text{implies}\quad\H\al\subseteq\H\beta.
\end{equation}

\begin{proposition}\label{prop:F-al}
\emph{\cite{binom}}\ 
For natural $\al$, one has $f\in\H\al$ if and only if $f$ has finite derivatives $f^{(0)}:=f,f^{(1)}:=f',\dots,f^{(\al-1)}$ on $\R$ such that $f^{(\al-1)}$ is convex on $\R$ and $f^{(j)}(-\infty+)=0$ for $j=0,1,\dots,\al-1$. 
\end{proposition}

It follows from Proposition \ref{prop:F-al} that, for every $t\in\R$, every $\beta\ge\al$, and every $\la>0$, the functions $u\mapsto(u-t)_+^\beta$ and $u\mapsto e^{\la(u-t)}$
belong to $\H\al$.  

Eaton (1970) \cite{eaton1} proved the Khinchin-Whittle-Haagerup inequality \eqref{eq:khin} for a class of moment functions, which essentially coincides with the class $\H3$; see \cite[Proposition~A.1]{pin94}.
Based on asymptotics, numerics, and a certain related inequality, 
Eaton (1974) \cite{eaton2} conjectured that the mentioned moment comparison inequality of his implies that 
$$\P\lp\vp_1a_1+\dots+\vp_n a_n\ge x\rp\le \frac{2e^3}9\,\frac1{x\sqrt{2\pi}}e^{-x^2/2}\quad \forall x>\sqrt2.
$$
Pinelis (1994) \cite{pin94} proved the following improvement of this conjecture: 
\begin{equation}\label{eq:pin94}
\P\lp\vp_1a_1+\dots+\vp_n a_n\ge x\rp\le\frac{2e^3}9\,\P(Z\ge x)\quad
\forall x\in\R. 
\end{equation}

It was realized in Pinelis (1998) \cite{pin98} that the reason why it is possible to extract tail comparison inequality \eqref{eq:pin94} from the Khinchin-Eaton moment comparison inequality \eqref{eq:khin} for $f\in\H3$ is that the tail function $x\mapsto\P(Z\ge x)$ is log-concave.
This realization resulted in a general device, which allows one to extract an optimal tail comparison inequality from an appropriate moment comparison inequality. The following is a special case of Theorem~4 of Pinelis (1999) \cite{pin99}; see also Theorem~3.11 of Pinelis (1998) \cite{pin98}.
\begin{theorem}
\label{th:comparison} 
Suppose that $0\le\beta\le\al$, $\xi$ and $\eta$ are real-valued r.v.'s, and the tail function $u\mapsto\P(\eta\ge u)$ is log-concave on $\R$. Then the comparison inequality 
\begin{equation}\label{eq:comp-al}
\E f(\xi)\le\E f(\eta)\quad\text{for all }f\in\H\al
\end{equation}
implies
\begin{equation}\label{eq:comp-beta}
\E f(\xi)\le c_{\al,\beta}\,\E f(\eta)\quad\text{for all }f\in\H\beta
\end{equation}
and, in particular, for all real $x$,
\begin{align}
\P(\xi\ge x) &\le \inf_{f\in\H\al}\,\frac{\E f(\eta)}{f(x)} \label{eq:comp-prob1} \\
&= B_{\opt}(x):=\inf_{t\in(-\infty,x)}\,\frac{\E(\eta-t)_+^\al}{(x-t)^\al} 
\label{eq:comp-prob2} \\
&\le \min\lp c_{\al,0}\,\P(\eta\ge x),\; \inf_{h>0}\,e^{-hx}\,\E e^{h\eta} \rp,
\label{eq:comp-prob3} 
\end{align}
where 
\begin{equation}\label{eq:c(al,beta)}
c_{\al,\beta}:=\frac{\Gamma(\al+1)(e/\al)^\al}{\Gamma(\beta+1)(e/\beta)^\beta}
\end{equation}
for $\beta>0$; $c_{\al,0}:=\Gamma(\al+1)(e/\al)^\al$.
Moreover, the constant $c_{\al,\beta}$ is the best possible in \eqref{eq:comp-beta} and \eqref{eq:comp-prob3}. 
\end{theorem}

A similar result for the case when $\al=1$ and $\beta=0$ 
is contained in the book by Shorack and Wellner (1986) \cite{shorack-wellner}, pages 797--799. 

Note that $c_{\al,0}\sim\sqrt{2\pi\al}$ and $c_{\al,\beta}\sim\sqrt{\al/\beta}$ as $\al,\beta\to\infty$.

\begin{remark}
\label{comparison-remark}
As follows from \cite[Remark 3.13]{pin98}, a useful point is that the requirement of the log-concavity of the tail function $q(u):=\P(\eta\ge u)$ in Theorem \ref{th:comparison}
can be relaxed by replacing $q$ with any (e.g., the least) log-concave majorant of $q$. 
However, then the optimality of $c_{\al,\beta}$ is not guaranteed. 
\end{remark}

Detailed studies of various cases and aspects of the optimal bound $B_{\opt}(x)$ in \eqref{eq:comp-prob2} were presented in \cite{dufour,pin98,be-64pp}.

Note that 
$c_{3,0}=2e^3/9$, 
which is the constant factor in \eqref{eq:pin94}.
Bobkov, G\"{o}tze, and Houdr\'{e} (BGH) (2001) \cite{houdre} obtained a simpler proof of inequality \eqref{eq:pin94},
but with a constant factor $\approx12.0$ in place of $2e^3/9\approx4.46$. 
In \cite{pin-BGH} the BGH method was modified to obtain a version of \eqref{eq:pin94} with a constant factor $\approx3.22$, which is $\approx1.01$ times the least possible constant factor in \eqref{eq:pin94}. 
Edelman \cite{edel} proposed inequality
$\P(S_n\ge x) \le \P\lp Z\ge x-1.5/x\rp$ for all $x>0$, but his proof appears to have a gap. A more precise upper bound, with $\ln c_{3,0}=1.495\dots$ in place of $1.5$, was recently shown \cite{pin-ed} to be a rather easy corollary of \eqref{eq:pin94}.

\begin{remark}\label{rem:hilbert}
One also has two kinds of multi-dimensional analogues of \eqref{eq:khin} and \eqref{eq:pin94}. One kind is represented by \cite[Theorems~2.3 and 2.4]{pin94}. The other kind is based on the dimensionality reduction device given in \cite{dim-reduct}.  
Indeed, Remarks in \cite{eaton2} imply (cf. the proof of Lemma~3.2 in \cite{pin94}) that, for any even function $f$ in class $\FF3$ (which contains $\H3$ and is defined by \eqref{eq:FF3} below), the function $[0,\infty)\ni u\mapsto f(\sqrt u)$ is convex.
Therefore, by \cite[Theorem~2.1]{dim-reduct}, 
\begin{equation}\label{eq:khin-hilbert}
\E f\lp \|\vp_1\xx_1+\dots+\vp_n\xx_n\| \rp\le\E f(|Z|)\quad
\text{for all even $f$ in $\FF3$},
\end{equation}
where 
$\xx_1,\dots,\xx_n$ are any non-random vectors in a Hilbert space $(H,\|\cdot\|)$ such that $\|\xx_1\|^2+\dots+\|\xx_n\|^2=1$. 
It follows that 
\eqref{eq:khin-hilbert} holds for all functions $f$ given by $f(x):=\E g(\vp x)$ $(x\in\R)$, where $g\in\H3$ and $\vp$ is a Rademacher r.v. Hence,
by Theorem~\ref{th:comparison},
\begin{equation}\label{eq:pin94-hilbert}
\P\lp  \|\vp_1\xx_1+\dots+\vp_n\xx_n\| \ge x \rp\le\frac{2e^3}9\,\P(|Z|\ge x)\quad
\forall x\in\R, 
\end{equation}
More generally, in view of a result by Hunt \cite{hunt}, inequalities 
\eqref{eq:khin}, \eqref{eq:pin94}, \eqref{eq:khin-hilbert}, and \eqref{eq:pin94-hilbert} hold if $\vp_1,\dots,\vp_n$ are replaced there by any independent zero-mean r.v.'s $\eta_1,\dots,\eta_n$ such that $|\eta_i|\le1$ almost surely (a.s.) for all $i$.  
\end{remark}

Pinelis (1999) \cite{pin99} also obtained the ``discrete" improvement of \eqref{eq:pin94}: 
\begin{equation}\label{eq:pin99}
\P\lp\vp_1a_1+\dots+\vp_n a_n\ge x\rp 
\le\frac{2e^3}9\,
\P\lp\frac1{\sqrt n}(\vp_1+\dots+\vp_n)\ge x\rp 
\end{equation}
for all values $x$ of r.v. $\frac1{\sqrt n}(\vp_1+\dots+\vp_n)$.

Such results can be, and have been, extended in several different directions.
In what follows, let $(S_0,S_1,\dots)$ be a supermartingale relative to a nondecreasing sequence of $\sigma$-algebras $(H_{\le0},H_{\le1},\dots)$, with $S_0\le0$ a.s. and differences 
$$X_i:=S_i-S_{i-1},\quad i=1,2,\dots.$$ 

The following normal domination statement is one of the main results of \cite{normal}.

\begin{theorem}\label{th:bernoulli} 
\emph{\cite{normal}}\ 
Suppose that for every $i=1,2,\dots$ there exist $H_{\le(i-1)}$-measurable r.v.'s $A_{i-1}$ and $B_{i-1}$ and a positive real number $c_i$ such that 
\begin{gather}
-A_{i-1}\le X_i\le B_{i-1} \quad\text{ and } \label{eq:bern-cond1} \\
\tfrac12(A_{i-1}+B_{i-1})\le c_i \label{eq:bern-cond2}
\end{gather}
a.s. Then for all $f\in\H5$ and all $n=1,2,\dots$
\begin{equation}\label{eq:bernoulli}
\E f(S_n)\le\E f(s\sqrt{n}Z),
\end{equation}
where 
\begin{equation}\label{eq:s-range}
s:=s^{(1)}:=\sqrt{\frac{c_1^2+\dots+c_n^2}{n}}.
\end{equation}
\end{theorem}

Note that inequality \eqref{eq:bernoulli} for the smaller class of exponential functions in place of the class $\H5$ is due to Hoeffding \cite{hoeff}.

By virtue of Theorem \ref{th:comparison}, one has the following corollary. 

\begin{corollary}\label{cor:bern-beta}
\emph{\cite{normal}}\ 
Under the conditions of Theorem~\ref{th:bernoulli}, for all $\beta\in[0,5]$, all $f\in\H\beta$, and all $n=1,2,\dots$
\begin{equation}\label{eq:bern-beta}
\E f(S_n)\le c_{5,\beta}\,\E f(s\sqrt{n}Z).
\end{equation}
In particular, for all real $x$,
\begin{align}
\P(S_n\ge x) &\le \inf_{f\in\H5}\,\frac{\E f(s\sqrt{n}Z)}{f(x)} \label{eq:bern-prob1} \\
&= \inf_{t\in(-\infty,x)}\,\frac{\E(s\sqrt{n}Z-t)_+^5}{(x-t)^5} 
\label{eq:bern-prob2} \\
&\le \min\lp c_{5,0}\P(s\sqrt{n}Z\ge x), \inf_{h>0}\,e^{-hx}\,\E e^{hs\sqrt{n}Z} \rp 
\label{eq:bern-prob3} \\
&= \min\lp c_{5,0}\P\lp Z\ge\tfrac{x}{s\sqrt{n}}\rp, \exp\lp -\tfrac{x^2}{2ns^2} \rp\rp. 
\label{eq:bern-prob4} 
\end{align}
\end{corollary}

The upper bound \eqref{eq:bern-prob4} -- but with a constant factor greater than
$427$ 
in place of $c_{5,0}=5!(e/5)^5=5.699\dots$ was obtained in Bentkus (2001) \cite{bent-isr} for the case when $(S_i)$ is a martingale. 
(In this case, Bentkus was using direct methods, rather than a generalized moment comparison inequality such as \eqref{eq:bernoulli}.)
The large value, $427$, of the constant factor renders the bound in \cite{bent-isr} hardly usable in statistics. Indeed, the upper bound
$427\P\lp Z\ge\tfrac{x}{s\sqrt n}\rp$ improves the Hoeffding bound
$\exp\left(-\frac{x^2}{2ns^2}\right)$ 
only when $\frac x{s\sqrt n}>170$, in which case (in view of \eqref{eq:bern-prob4}) one has 
$\P(S_n\ge x)<c_{5,0}\P\lp Z\ge170\rp<10^{-6200}$.



As shown in \cite{normal}, Theorem~\ref{th:bernoulli} and Corollary~\ref{cor:bern-beta} are well suited in order to obtain the most precise presently known bounds for the measure concentration phenomenon in terms of separately-Lipschitz (or, equivalently, $\ell^1$-Lipschitz) functions on product spaces. 

Theorem~\ref{th:bernoulli} can be further improved, as follows.

\begin{theorem}\label{th:bernoulli-improved} 
\emph{\cite{normal}}\ 
Suppose that for every $i=1,2,\dots$ there exist a positive $H_{\le(i-1)}$-measurable r.v. $B_{i-1}$ and a positive real number $\hat s_i$ such that 
\begin{gather}
X_i\le B_{i-1} \quad \text{and} \label{eq:bern-improv-cond1} \\ \frac12\left(B_{i-1}+\frac{\Var_{i-1}X_i}{B_{i-1}}\right)\le \hat s_i \label{eq:bern-improv-cond2} 
\end{gather}
a.s. 
Here and elsewhere, we let $\Var_j$ stand for conditional variance given $H_{\le j}$. 
Then one has all the inequalities \eqref{eq:bernoulli} and \eqref{eq:bern-beta}--\eqref{eq:bern-prob4}, only with $s$ replaced by  
\begin{equation}\label{eq:hat-s}
\hat s:=\sqrt{\frac{\hat s_1^2+\dots+\hat s_n^2}n}.
\end{equation} 
\end{theorem}

The set of conditions \eqref{eq:bern-improv-cond1}--\eqref{eq:bern-improv-cond2} is equivalent to
$$X_i\le B_{i-1} \quad \text{and}\quad \si_*(B_{i-1}, \Var_{i-1}X_i) \le\hat s_i
$$
a.s., where
\begin{multline*}
\si_*(b_0,c^2):=
\frac12\,\inf_{b\ge b_0} \left(b+\frac{ c^2 }b \right)
=\min\left( c\vee b_0, \frac12\left(b_0+\frac{ c^2 }{b_0}\right)\right) \\
=\begin{cases}
c & \text{if } c\ge b_0, \\
\frac12\, \left(b_0+\frac{ c^2 }{b_0}\right) & \text{if } c<b_0, 
\end{cases}
\end{multline*}
for positive $c$ and $b_0$. 
This follows simply because the inequalities $X_i\le B_{i-1}$ and $b\ge B_{i-1}$ imply $X_i\le b$.  

Thus, in the case when $\Var_{i-1}X_i<B_{i-1}^2$ a.s., conditions \eqref{eq:bern-improv-cond1}--\eqref{eq:bern-improv-cond2} represent an improvement of condition $B_{i-1}^2\vee\,\Var_{i-1}X_i\le\hat s_i^2$ a.s., imposed in \cite{bent-jtp,bent-ap}. In a certain variety of cases, this improvement may be even more significant than the improvement in the constant factor from $427$ to $5.699\dots$ \emph{before} the probability sign. 

Moreover, it can be shown that the function $\si_*(\cdot,\cdot)$ of the pair \\
$(B_{i-1}, \Var_{i-1}X_i)$, which is in effect used in Theorem~\ref{th:bernoulli-improved} is nearly optimal as far as the normal domination is concerned. 

However, it can also be seen that even the best possible normal domination may be inadequate if the asymmetry of the random summands $X_i$ is significant or if $n$ is not large. 
In such a case, one may try to use binomial domination instead of normal, as in the following theorem, which is a straighforward corollary of results of \cite{bent-liet02} (or \cite{bent-ap}). 

\begin{theorem}\label{th:bent}
\emph{\cite{bent-liet02,bent-ap}}
Suppose that for every $i=1,2,\dots$ there exist non-random constants $b_i>0$ and $c_i>0$ 
such that 
\begin{align}
X_i&\le b_i \quad\text{and } \label{eq:bent-cond1} \\
\Var_{i-1}X_i&\le c_i^2 \label{eq:bent-cond2}
\end{align}
a.s. Then, for all $n=1,2,\dots$,
\begin{align}\label{eq:bent}
\E f(S_n)&\le\E f(T_n)\quad\forall f\in\H2,\quad\text{where} \\
T_n&:=Z_1+\dots+Z_n
\end{align}
and $Z_1,\dots,Z_n$ are i.i.d.\ r.v.'s such that each $Z_i$ takes on only two values, one of which is 
\begin{equation}\label{eq:d}
b:=\max_i b_i,
\end{equation}
and satisfies the conditions
\begin{gather}
\E Z_i=0\quad\text{and}\quad\Var Z_i=c^2,\quad\text{where} \notag\\
c:=\lp\dfrac1n\sum_{i=1}^nc_i^2\rp^{1/2};\label{eq:si}
\end{gather} 
$$\quad\text{that is,}\quad
\PP(Z_i=b)=\frac{c^2}{b^2+c^2}\quad\text{and}\quad
\PP\lp Z_i=-\frac{c^2}{b}\rp=\frac{b^2}{b^2+c^2}.$$ 
\end{theorem}

Based on this result, the tail comparison inequality
\begin{equation}\label{eq:interp}
\PP(S_n\ge y) \le c_{2,0} \PP^{\lin,\lc}(T_n\ge y+\tfrac h2)\quad\forall y\in\R
\end{equation}
was obtained in \cite{binom},
where $c_{2,0}=e^2/2=3.694\dots$ (in accordance with \eqref{eq:c(al,beta)}), 
$h:=b+c^2/b$, and the function $y\mapsto\PP^{\lin,\lc}(T_n\ge y)$ is the least log-concave majorant of the linear interpolation of the tail function $y\mapsto\PP(T_n\ge y)$ over the lattice of all points of the form $nb+kh$ ($k\in\Z$). Tail comparison inequality \eqref{eq:interp} is a substantial improvement of the corresponding inequality in \cite[Theorem~1]{bent-liet02} and \cite[Theorem~1.1]{bent-ap}.

Yet, while the ``variance"-averaging given by \eqref{eq:si} is nice, the extreme kind of upper-bound-averaging \eqref{eq:d} seems very crude. 

In this paper, another approach to the problem of asymmetry is presented. Here we provide binomial upper bounds on generalized moments and tail probabilities for $S_n$ assuming that certain \emph{indices of asymmetry} of the $X_i$'s (rather than the $X_i$'s themselves) are \emph{uniformly} bounded from above. 
This assumption of bounded asymmetry (in contrast with the uniform boundedness) of the $X_i$'s is natural in certain statistical applications; see Subsection~\ref{t-stats}.

The results of \cite{bent-symm} for the symmetric case can be similarly complemented, with condition \eqref{eq:d} of the uniform boundedness of the ranges now replaced by a condition of uniform boundedness of the kurtoses of the $X_i$'s; see Remark~\ref{rem:symm}.


\section{Statements of basic results and discussion}\label{results}

Let $\bs(p)$ denote the standardized Bernoulli distribution with parameter $p$: for a r.v.\ $\BS$ let, by definition,
$$\BS\sim\bs(p) \iff \P\lp\BS=\sqrt{\tfrac qp}\,\rp=p=1-\P\lp\BS=-\sqrt{\tfrac pq}\,\rp,$$
where
$$q:=1-p\quad\text{and}\quad 0<p<1;$$
thus, $\bs(p)$ is a two-point zero-mean unit-variance distribution. 
In particular, $\bs(\frac12)$ is the distribution of a Rademacher r.v.\ $\vp$, with $\P(\vp=\pm1)=\frac12$. 

Let $\C^2$ denote the class of all twice continuously differentiable functions 
$f\colon\R\to\R$.
Consider the following class of functions: 
\begin{equation} \label{eq:F3}
\F3:=\{f\in\C^2\colon \text{$f$ and $f''$ are nondecreasing and convex}\}.
\end{equation} 
An equivalent definition would be given by the formula
$$\F3=\{f\in\C^2\colon \text{$f$, $f'$, $f''$, and $f'''$ are nondecreasing}\},$$
where $f'''$ denotes the right derivative of the convex function $f''$. 

For example, functions $x\mapsto a+b\,x+c\,(x-t)_+^\al$ and 
$x\mapsto a+b\,x+c\,e^{\la x}$ belong to $\F3$ for all $a\in\R$, $b\ge0$, $c\ge0$, $t\in\R$, $\al\ge3$, and $\la\ge0$. 

\begin{remark} \label{rem:Jensen's}
If a function $f\colon\R\to\R$ is convex and a r.v.\ $X$ has a finite expectation, then, by Jensen's inequality, $\E f(X)$ always exists in $(-\infty,\infty]$. This remark will be used in this paper (sometimes tacitly) for functions $f$ in the class $\F3$, as well as for other convex functions. 
\end{remark} 

Throughout the paper, unless indicated otherwise, the following notation/\- assumptions will be used: 
\begin{equation} \label{eq:restrs}
m\in[1,\infty),\quad p\in(0,1),\quad q=1-p,\quad\text{and}\quad
\BS_1,\dots,\BS_n\iid\bs(p).
\end{equation}

Introduce also
\begin{equation} \label{eq:m_*}
m_*(p):=
\begin{cases}
\dfrac{1 + p + 2\,p^2}
   {2{\left( {\sqrt{p - p^2}} + 
       2\,p^2 \right) }} \quad & \text{if}\quad 0<p\le\frac{1}{2}, \\
1 \quad & \text{if}\quad \frac{1}{2}\le p<1.     
\end{cases}
\end{equation}
Later it will be clear that $m_*(p)$ increases from $1$ to $\infty$ as $p$ decreases from $\frac12$ to $0$ (see the proof of Lemma~\ref{lem:p*,m*}). 

Of the main results of this paper, the following one is perhaps the easiest to state (but not to prove).

\begin{theorem}\label{th:bern}
For any real number 
$$m\ge m_*(p),$$
all $f\in\F3$, 
all natural $n$, and all nonnegative numbers 
$c_1,\dots,c_n$, one has
\begin{gather}
\E f(c_1\BS_1+\dots+c_n\BS_n)\le\E f\lp s^{(m)}\cdot(\BS_1+\dots+\BS_n)\rp,
\label{eq:Ef bern}\\
\quad\text{where}\quad s^{(m)}:=\lp\frac1n\,\sum_{i=1}^n c_i^{2m}\rp^\frac1{2m}.
\label{eq:s(m)}
\end{gather}
Moreover, the lower bound $m_*(p)$ on $m$ is exact for each $p\in(0,1)$.
\end{theorem}

The proofs are deferred to Section~\ref{proofs}.

\begin{remark}\label{rem:m>1} 
The general restriction $m\ge1$ in \eqref{eq:restrs} is quite natural. Indeed, if inequality \eqref{eq:Ef bern} held for some $m\in(0,1)$ then, taking $c_1=1$, $c_2=\dots=c_n=0$, and letting $n\to\infty$, one would have, by the central limit theorem, the inequality $\E f(\BS_1)\le f(0)$ for all $f\in\F3$, which is false even for $f(x)\equiv e^x$ or $f(x)\equiv x_+^3$.
\end{remark}

Here is a generalization of Theorem~\ref{th:bern}: 

\begin{theorem}\label{th:1} 
Let $X_1,\dots,X_n$ be independent r.v.'s such that for every $i\in\{1,\dots,n\}$
$$\E X_i\le0\quad\text{and}\quad -a_i\le X_i\le b_i \ 
\text{a.s.},$$
where $a_i$ and $b_i$ are positive numbers such that
\begin{equation}\label{eq:bi/ai}
\frac{b_i}{a_i}\le\frac qp.
\end{equation}
Then, for any real number 
$m\ge m_*(p)$
and all $f\in\F3$, 
one has the inequality
\begin{gather}
\E f(X_1+\dots+X_n)\le\E f\lp s^{(m)}\cdot(\BS_1+\dots+\BS_n)\rp,
\quad\text{where}
\label{eq:Ef main} \\
s^{(m)}:=\lp\frac1n\sum_{i=1}^n (a_i b_i)^m\rp^\frac1{2m}.\notag
\end{gather}
 
Moreover, the lower bound $m_*(p)$ on $m$ is exact for each $p\in(0,1)$.
\end{theorem}
 
Condition \eqref{eq:bi/ai} may be referred to as a bounded-asymmetry coundition.

Theorem~\ref{th:1} can be easily extended to (super)martingales. 

\begin{theorem}\label{th:superm} 
Suppose that for every $i\in\{1,\dots,n\}$ one has \eqref{eq:bern-cond1} with positive $A_{i-1}$ and $B_{i-1}$ such that
\begin{align}
\sqrt{A_{i-1}\,B_{i-1}} & \le c_i\quad\text{and} \label{eq:AiBi} \\
\frac{B_{i-1}}{A_{i-1}} & \le\frac qp \label{eq:Bi/Ai} 
\end{align}
a.s., where $c_i$ is a non-random number. 
Then, for any real number 
\begin{equation}\label{eq:m>m*(p)}
m\ge m_*(p)
\end{equation}
and $f\in\F3$,
one has the inequality
\begin{gather}
\E f(S_n)\le\E f\lp s^{(m)}\cdot(\BS_1+\dots+\BS_n)\rp, 
\label{eq:Ef superm} 
\end{gather}
where $s^{(m)}$ is defined by \eqref{eq:s(m)}.
Moreover, the lower bound $m_*(p)$ on $m$ is exact for each $p\in(0,1)$.
\end{theorem}

One should compare \eqref{eq:AiBi} and \eqref{eq:s(m)} with \eqref{eq:bern-cond2} and \eqref{eq:s-range}. If $X_{a,b}$ stands for a zero-mean r.v.\ taking on values in the set $\{-a,b\}$ for some positive $a$ and $b$, then
obviously the half-range $\frac12(a+b)$ of $X_{a,b}$ is no less than its standard deviation $\sqrt{ab}$. That is, \eqref{eq:bern-cond2} is more restrictive than \eqref{eq:AiBi}. On the other hand, one has the inequality $s^{(m)}\ge s^{(1)}$ for $m\ge1$. Moreover, the greater the uniform bound $\frac qp$ on asymmetry in \eqref{eq:Bi/Ai} is, the greater $m$ must be according to \eqref{eq:m>m*(p)} and hence the more pronounced the inequality $s^{(m)}\ge s^{(1)}$ will be.
Yet, it will be demonstrated elsewhere that, overall, \eqref{eq:AiBi} and \eqref{eq:s(m)} work better in certain important statistical applications than \eqref{eq:bern-cond2} and \eqref{eq:s-range}.
Note also that one can choose the ``ideal'' value $m=1$ whenever the asymmetry index $\frac qp$ does not exceed $1$, that is, whenever the $X_i$'s are not skewed to the right. 

Recall the definition of the
{\em Schur majorizarion\/}: for $\mathbf{a:=}\left( a_{1},\ldots,a_{n}\right) $ and $%
\mathbf{b:=}\left( b_{1},\ldots,b_{n}\right) $ in $\mathbb{R}^{n}$, $\mathbf{%
a\succeq b}$ means that $a_{1}+\cdots+a_{n}=b_{1}+\cdots+b_{n}$ and $a_{%
\left[ 1\right] }+\cdots+a_{\left[ j\right] }\geq b_{\left[ 1\right]
}+\cdots+b_{\left[ j\right] }$ for all $j\in\left\{ 1,\ldots,n\right\} $,
where $a_{\left[ 1\right] }\geq\cdots\geq a_{\left[ n\right] }$ are the
ordered numbers $a_{1},\ldots,a_{n}$, from the largest to the smallest.  
Recall also that a function 
$\mathcal{Q}\colon[0,\infty)^n\rightarrow\mathbb{R}$ is referred to as {\em Schur-concave} if it reverses the
Schur majorization: for any $\mathbf{a}$ and $\mathbf{b}$ in $[0,\infty)^n$ such that 
$\mathbf{a}\succeq\mathbf{b}$, one has 
$\mathcal{Q}\left( \mathbf{a}\right)
\leq\mathcal{Q}\left( \mathbf{b}\right)$. 

Theorems~\ref{th:bern}, \ref{th:1}, and \ref{th:superm} are contained in 

\begin{theorem}\label{th:main} 
The following statements are equivalent to one another.
\emph{
\begin{enumerate}[(I)]
	\item \label{it:m>m_*(p)} $m\ge m_*(p)$.
	\item \label{it:Bern}  \emph{For all $f\in\F3$, all natural $n\ge2$, and all nonnegative numbers 
$c_1,\dots,c_n$, one has \eqref{eq:Ef bern}.}
	\item \label{it:Bern,n=2} \emph{The same as item \eqref{it:Bern}, but only for $n=2$.} 
	\item \label{it:Bern-Sch} \emph{For every natural $n\ge2$ and every function 
$f\in\F3$, the function 
\begin{equation} \label{eq:schur} 
[0,\infty)^n\ni(a_1,\dots,a_n)\longmapsto
\E f(a_1^{1/(2m)}\BS_1+\dots+a_n^{1/(2m)}\BS_n)
\end{equation}
is Schur-concave.}
	\item \label{it:Bern-Sch,n=2} \emph{The same as item \eqref{it:Bern-Sch}, but only for $n=2$.}  
	\item \label{it:bounded} 
\emph{Let the $X_i$'s, $a_i$'s, $b_i$'s, and 
$s^{(m)}$ be as in the statement of Theorem~\ref{th:1}. Then one has \eqref{eq:Ef main} for all natural $n$.} 
	\item \label{it:superm} 
\emph{Let the $X_i$'s, $A_{i-1}$'s, $B_{i-1}$'s, 
$c_i$'s, and $s^{(m)}$ be as in the statement of Theorem~\ref{th:superm}. Then one has 
\eqref{eq:Ef superm} for all natural $n$. }
\end{enumerate}
}
\end{theorem}

The special case of statement \eqref{it:Bern-Sch} of Theorem~\ref{th:main} with $p=\frac12$ and $m=1$ is essentially the mentioned result due to Whittle~\cite{whittle} and Eaton~\cite{eaton1}.

From the ``right-tail" Theorem~\ref{th:main}, one can deduce its left-tail and two-tail analogues. Appropriate left-tail and two-tail counterparts of $\F3$ are the following classes of functions:
\begin{align} 
\Fminus3 & :=\{f\in\C^2\colon \text{$f$ and $f''$ are nonincreasing and convex}\}
\label{eq:Fminus3} \\
& =\{f\colon \exists g\in\F3\ \forall x\in\R\ f(x)=g(-x)\}
\quad\text{and} \notag\\
\FF3 & :=\{f\in\C^2\colon \text{$f$ and $f''$ are convex}\} \label{eq:FF3}.
\end{align} 

\begin{theorem}\label{th:main-lefttail} 
Theorem~\ref{th:main} holds with $\Fminus3$ in place of $\F3$ if the restrictions 
\emph{(i)} $m\ge m_*(p)$, 
\emph{(ii)} $\E X_i\le0$, 
\emph{(iii)} $\dfrac{b_i}{a_i}\le \dfrac qp$, 
\emph{(iv)} $(S_0,\dots,S_n)$ is a supermartingale with $S_0\le0$, and 
\emph{(v)} $\dfrac{B_{i-1}}{A_{i-1}}\le \dfrac qp$ in Theorem~\ref{th:main} are replaced, respectively, with the following:
\emph{
\begin{enumerate}[(i)]
\item $m\ge m_*(q)$, 
\item $\E X_i\ge0$, 
\item $\dfrac{a_i}{b_i}\le \dfrac pq$,
\item \emph{$(S_0,\dots,S_n)$ is a submartingale with $S_0\ge0$, and}
\item 
$\dfrac{A_{i-1}}{B_{i-1}} \le \dfrac pq$. 
\end{enumerate}
}
\end{theorem}

This ``left-tail" analogue is a trivial corollary of Theorem~\ref{th:main}. 

The ``two-tail'' analogue of Theorem~\ref{th:main} is more difficult to prove. 
It relies in part on Proposition~\ref{prop:f3,ff3} below, preceded by the following definition.

\begin{definition} \label{def:convergence}
Let us say that a sequence of functions $(f_n)$ in $\C^2$ converges to a function $f$ in $\C^2$ and write $f_n\to f$ (as $n\to\infty$) if $f_n(x)\uparrow f(x)$ and 
$f''_n(x)\to f''(x)$ for all real $x$. 
(This stronger notion of convergence will make it easier to verify the convergence of relevant expected values; also, it naturally provides for the relevant classes of functions to be closed.) 

For any subset $\A$ of $\C^2$, its {\em closure} -- denoted here by $\cl\A$ -- will be understood here simply as the set of the limits of all sequences in $\A$ that are convergent in $\C^2$. Obviously, $\cl\A\supseteq\A$, for all $\A\subseteq\C^2$.
\end{definition}

Obviously, the ``two-tail'' class $\FF3$ contains both ``one-tail" classes $\F3$ and 
$\Fminus3$.
The more informative relation of $\FF3$ to $\F3$ and 
$\Fminus3$ (on which the proof of Theorem~\ref{th:main-2tail} below is partly based) is given by 

\begin{proposition} \label{prop:f3,ff3}
One has $\FF3=\cl\GG3$, where
\begin{multline} \label{eq:GG3}
\GG3:=\{f\in\C^2\colon \exists c\ge0\ \exists f_+\in\F3\ \exists f_-\in\Fminus3\ \forall x\in\R\ \\
f(x)=c\,x^2/2+f_+(x)+f_-(x)\}.
\end{multline}
However, $\FF3\ne\GG3$.
\end{proposition}

For example, functions $x\mapsto a+b\,x+c\,x^2+d\,|x-t|^\al$, 
$x\mapsto\cosh{\la x}$, $x\mapsto e^{\la x}$, $x\mapsto(x-t)_+^\al$, and 
$x\mapsto(t-x)_+^\al$ belong to $\FF3$ for all $a\in\R$, $b\in\R$, $c\ge0$, $d\ge0$, 
$t\in\R$, $\al\ge3$, and $\la\in\R$.

Note also that the classes $\F3$, $\Fminus3$, and $\FF3$ are convex cones; that is, any linear combination with nonnegative coefficients of functions belonging to any one of these classes belongs to the same class.

\begin{theorem}\label{th:main-2tail} 
Theorem~\ref{th:main} holds with $\FF3$ in place of $\F3$ if the restrictions 
\emph{(i)} $m\ge m_*(p)$, 
\emph{(ii)} $\E X_i\le0$, 
\emph{(iii)} $\frac{b_i}{a_i}\le \frac qp$, 
\emph{(iv)} $(S_0,\dots,S_n)$ is a supermartingale with $S_0\le0$, and 
\emph{(v)} $\frac{B_{i-1}}{A_{i-1}}\le \frac qp$ in Theorem~\ref{th:main} are replaced, respectively, with the following stronger restrictions:
\emph{
\begin{enumerate}[(i)]
\item \emph{$m\ge m_*(p) $ and $p\le\frac12$,} 
\item $\E X_i=0$, 
\item $\max\lp \dfrac{b_i}{a_i},\dfrac{a_i}{b_i}\rp\le \dfrac qp$,
\item \emph{$(S_0,\dots,S_n)$ is a martingale with $S_0=0$, and}
\item 
$\max\lp \dfrac{B_{i-1}}{A_{i-1}},\dfrac{A_{i-1}}{B_{i-1}}\rp \le \dfrac qp$. 
\end{enumerate}
}
\end{theorem}

\section{Applications}\label{appls} 

\subsection{Bounds on even richer classes of generalized moments (including tail probabilities), maximal inequalities, and some further extensions}\label{beta,tails} 

Using Theorem~\ref{th:superm} and Remark~\ref{comparison-remark} (and also recalling definition \eqref{eq:F3} and Proposition~\ref{prop:F-al}), one immediately obtains the following corollary, which may be compared with Corollary~2.2 in \cite{normal}. 

\begin{corollary}\label{cor:beta,prob}
Suppose that conditions \eqref{eq:bern-cond1}, \eqref{eq:AiBi}, \eqref{eq:Bi/Ai}, and \eqref{eq:m>m*(p)} hold. 
Then for all $\beta\in[0,3]$, all $f\in\H\beta$, and all $n=1,2,\dots$
\begin{equation}\label{eq:bern-3-beta}
\E f(S_n)\le c_{3,\beta}\E f(T_n),
\end{equation}
where
$$T_n:=s^{(m)}\cdot(\BS_1+\dots+\BS_n).$$
In particular, for all real $x$,
\begin{align}
\P(S_n\ge x) &\le \inf_{f\in\H3}\,\frac{\E f(T_n)}{f(x)} \label{eq:prob1} \\
&= \inf_{t\in(-\infty,x)}\,\frac{\E(T_n-t)_+^3}{(x-t)^3} 
\label{eq:prob2} \\
&\le \min\lp c_{3,0}\,\P^\lc(T_n\ge x), \inf_{\la>0}\,e^{-\la x}\,\E e^{\la T_n} \rp,
\label{eq:prob3} \\
&= \min\lp c_{3,0}\,\P^\lc(T_n\ge x), e^{-nH}\rp,
\label{eq:prob4} 
\end{align}
where $x\mapsto\P^\lc(T_n\ge x)$ is the least log-concave majorant of the function $x\mapsto\P(T_n\ge x)$ on $\R$; 
$H:=(p+y)\ln\frac{p+y}p+(q-y)\ln\frac{q-y}q$ if 
$0\le y:=\frac xn\,\frac{\sqrt{pq}}{s^{(m)}}<q$, $H:=-\ln p$ if $y=q$, $H:=\infty$ if $y>q$, and $H:=0$ if $y<0$. 
\end{corollary} 

Note that $\P^\lc(T_n\ge x)=\P(T_n\ge x)$ for all $x$ in the lattice $$L:=\{nb+kh\colon k\in\Z\}$$
generated by the support of the distribution of $T_n$, where $b:=s^{(m)}\sqrt{\frac qp}$ and $h:=s^{(m)}/\sqrt{pq}$.

The bound $e^{-nH}$ in \eqref{eq:prob4} is largely due to Hoeffding~\cite{hoeff}.

Using also results of \cite{binom}, one has the following.
\begin{corollary}\label{cor:interp}
Under the conditions of Corollary~\ref{cor:beta,prob},
\begin{equation}\label{eq:prob5} 
\P(S_n\ge x) \le c_{3,0}\,\P^{\lin,\lc}(T_n\ge x+\tfrac h2)\quad\forall x\in\R,
\end{equation}
where 
$x\mapsto\P^{\lin,\lc}(T_n\ge x)$ is the least log-concave majorant of the linear interpolation of the tail function $x\mapsto\P(T_n\ge x)$ over the lattice $L$.
\end{corollary} 

The upper bound in \eqref{eq:prob5} usually works better than that in \eqref{eq:prob4} in statistical practice. An explicit formula for 
$\P^{\lin,\lc}(T_n\ge x+\tfrac h2)$ is given in \cite{binom}. 

\begin{corollary}\label{cor:max}
In view of results of \cite{normal}, one can replace $S_n$ in the left-hand side of inequalities \eqref{eq:prob1} and \eqref{eq:prob5} by 
$$M_n:=\max_{0\le k\le n}S_k.$$
Similarly, inequality \eqref{eq:bern-3-beta} holds for all $\beta\in(0,3)$ with $M_n$ in place of $S_n$ if $c_{3,\beta}$ is replaced there by 
$\dfrac{k_{1;3,\beta}}{k_{3,\beta}}\,c_{3,\beta}$, where, for $\beta\in(0,\al)$,
\begin{equation*}\label{eq:k1,k}
k_{1;\al,\beta}:=
\sup_{\si>0} \si^{-\beta(\al-1)}
\left(\int_0^\si \frac{\beta s^{\beta-1}\,ds}{1+s}\right)^\al
\quad\text{and}\quad
k_{\al,\beta}:=\frac{\beta^\beta(\al-\beta)^{\al-\beta}}{\al^\al}.
\end{equation*}
\end{corollary}

\begin{corollary}\label{cor:bent}
The set of conditions \eqref{eq:bern-cond1}, \eqref{eq:AiBi}, and \eqref{eq:Bi/Ai} in Corollaries~\ref{cor:beta,prob}, \ref{cor:interp}, and \ref{cor:max} can be replaced by the set of conditions \eqref{eq:bent-cond1}, \eqref{eq:bent-cond2}, and
\begin{equation}\label{eq:bent-asymm}
\dfrac{b_i^2}{c_i^2}\le\dfrac qp\quad\forall i=1,\dots,n.
\end{equation}
\end{corollary}

In fact, the two sets of conditions mentioned in Corollary~\ref{cor:bent} are equivalent to each other in a certain sense; see e.g. Remark~2.4 in \cite{normal} and the proof of Theorem~2.3 therein.

Note that the special case of Corollary~\ref{cor:bent} with $p=\frac12$ (so that one may take $m=1$) contains, among other things, Theorem~1.3 of \cite{bent-ap}, which states that, if $(S_i)$ is a martingale with $S_0=0$ satisfying conditions \eqref{eq:bent-cond1} with 
$$b_i=c_i\quad\forall i=1,\dots,n$$
and \eqref{eq:bent-cond2}, then $\forall x\in\R$
\begin{equation}\label{eq:bent-Th1.3}
\P(S_n\ge x)\le c_{3,0}\P^\lc(T_n\ge x)\quad\text{with}\quad
T_n=s^{(1)}\cdot(\vp_1+\dots+\vp_n);
\end{equation}
by the central limit theorem, this inequality implies
\begin{equation}\label{eq:bent-Th1.3,normal}
\P(S_n\ge x)\le c_{3,0}\P(s^{(1)}\sqrt{n}Z\ge x).
\end{equation}
Obviously, inequalities~\eqref{eq:bent-Th1.3,normal} and \eqref{eq:bent-Th1.3} are extensions of \eqref{eq:pin94} and \eqref{eq:pin99}, respectively. A version of inequality \eqref{eq:bent-Th1.3,normal}, 
with the larger constant factor $1/\P(Z>\sqrt3)=24.01\dots$ in place of $c_{3,0}=2e^3/9=4.46\dots$, appeared earlier in \cite{bent-jtp}. The improvement in the constant factor achieved in \cite{bent-ap}, as compared with \cite{bent-jtp}, is due to replacing the direct method used in the earlier paper with the method based on Theorem~\ref{th:comparison} and Remark~\ref{comparison-remark}, which allows one to extract optimal tail comparison inequalities from comparison of generalized moments. 

More generally, the generalized moments $\E f(T_n)$ in the above upper bounds,
where $T_n:=s^{(m)}\cdot(\BS_1+\dots+\BS_n)$, can be replaced by $\E f(s^{(1)}\sqrt{n}Z)$ provided that $p\ge\frac12$.

That $(S_0,S_1,\dots)$ is allowed to be a supermartingale (rather than only a martingale) makes it convenient to use the simple but powerful truncation tool; cf. the discussion at the end of Section~2 in \cite{normal}.

\begin{remark}\label{rem:exp}
In the above results, the exact lower bound $m_*(p)$ on $m$ can be replaced by a substantially smaller (for $p\in(0,\frac12)$)
exact lower bound $m_{\mathsf{exp}}(p)$ on $m$ if 
the class $\F3$ is replaced by the substantially smaller class 
$$\F\infty:=\{f\in\C^\infty\colon f^{(j)}\ge0\ \forall j\in\{0,1,\dots\}\}$$
of completely monotone functions, where $\C^\infty$ is the class of all infinitely differentiable functions $f\colon\R\to\R$ with derivatives $f^{(0)}:=f, f^{(1)}:=f',\dots$. By Bernstein's theorem on completely monotone functions (see, e.g., \cite{choquet} or \cite{phelps}), 
$$\F\infty=\{f\colon f(x)=\smallint_{[0,\infty)}e^{kx}\mu(dk)\ \forall x\in\R\},$$
where $\mu\ge0$ is a Borel measure such that the integral $\int_{[0,\infty)}e^{kx}\mu(dk)$ is finite $\forall x\in\R$.
On the other hand, by \eqref{eq:F-al-beta} and Proposition~\ref{prop:F-al}, 
$$\F\infty=\{h+c\colon h\in\bigcap_{\al>0}\H\al,\ c\ge0\}.$$

By Remark~\ref{rem:m>1}, the exact lower bound $m_{\mathsf{exp}}(p)$ cannot be less than $1$, so that $m_{\mathsf{exp}}(p)=m_*(p)=1$ for all $p\in[\frac12,1)$. However, for $p\in(0,\frac12)$, $m_{\mathsf{exp}}(p)$ is substantially smaller than $m_*(p)$. In particular, when $p\downarrow0$, one has $m_{\mathsf{exp}}(p)\sim\frac12\ln\frac1p$, while $m_*(p)\sim\frac1{2\sqrt p}$. 
For $p\in(0,\frac12)$, the exact lower bound $m_{\mathsf{exp}}(p)$ can be described by parametric equations 
\begin{equation}\label{eq:param}
\left\{
\begin{aligned}
m_{\mathsf{exp}}(p)&=\tilde m(k):= \frac{(e^k+1) k}{2 (e^k-1)},
\\
p&=\tilde p(k):=\frac{e^k-1-k}{(e^k-1) (1+k+(k-1)e^k)}
,
\end{aligned}
\right.
\end{equation}
with $k\in(0,\infty)$ as the parameter. One can see that $\tilde p(k)$ decreases from $\frac12$ to $0$ as $k$ increases from $0$ to $\infty$. 
An explicit upper bound on $m_{\mathsf{exp}}(p)$ is given by the inequality
\begin{equation}\label{eq:m-exp-up}
m_{\mathsf{exp}}(p)<m_{\mathsf{exp,\,up}}(p):=
\frac{1-2p-\ln(2p)}{2(1-2p)}\quad\forall p\in(0,\tfrac12),
\end{equation}
so that one has $m\ge m_{\mathsf{exp}}(p)$ provided that $m\ge m_{\mathsf{exp,\,up}}(p)$. This simple upper bound on $m_{\mathsf{exp}}(p)$ is rather good: $m_{\mathsf{exp,\,up}}(p)\sim m_{\mathsf{exp}}(p)\sim\frac12\ln\frac1p$ as $p\downarrow0$, $m_{\mathsf{exp,\,up}}(p)\sim m_{\mathsf{exp}}(p)\sim1$ as $p\uparrow\frac12$, and $\max\limits_{0<p<1/2}\big(m_{\mathsf{exp,\,up}}(p)\,/ \,m_{\mathsf{exp}}(p)\big)=1.32\dots$, the maximum attained at $p=0.019\dots$.

On the other hand, the equation $\tilde p(k)=p$ has a unique solution  $k=:\tilde k(p)$ in $(0,\infty)$ for every $p\in(0,\frac12)$. Moreover, the Newton iterative scheme\\
$k_{j+1}:=k_j-F(k_j)/F'(k_j)$ for $j=0,1,\dots$, where 
$F(k):=F(p,k):=e^k-1-k \\ + (1+k-2e^k+e^{2 k}(1- k)) p  $, converges to $\tilde k(p)$ for every $p\in(0,\frac12)$ and every initial approximation $k_0>3(\frac1{\sqrt{2p}}-1)$. 
Thus, one has $m_{\mathsf{exp,\,up}}(p)=\tilde m(\tilde k(p))$ for all $p\in(0,\frac12)$. 
The following table illustrates the fact that $m_{\mathsf{exp}}(p)$ is substantially smaller than $m_*(p)$ for small $p$:
\begin{center}
	\begin{tabular}{l||l|l|l|l}
$p$				& $10^{-1}$	&	$10^{-2}$	&	$10^{-3}$	&	$10^{-4}$	\\
\hline 
$m_*(p)$	& $1.75$	& $5.06\dots$ & $15.83\dots$ & $50.00\dots$ \\	
$m_{\mathsf{exp}}(p)$	
					& $1.21\dots$	& $1.86\dots$ & $\ \;2.72\dots$ & $\ \;3.68\dots$
	\end{tabular}
\end{center}

In particular, it follows that if $p$ is small then the second upper bound in \eqref{eq:prob3}-\eqref{eq:prob4}, $e^{-nH}$, holds for significantly smaller values of $m$ than the first upper bound in \eqref{eq:prob3}-\eqref{eq:prob4}, $c_{3,0}\,\P^\lc(T_n\ge x)$, does. 
Thus, the exponential upper bound, $e^{-nH}$, may turn out to be smaller than the generally more precise upper bound, $c_{3,0}\,\P^\lc(T_n\ge x)$, even for large values of $x$ if the $X_i$'s differ very much from one another in distribution. 
Details on this remark, Remark~\ref{rem:exp}, will be presented elsewhere.
\end{remark}

\begin{remark}\label{rem:symm}
The stated results for the asymmetric case have ``symmetric" counterparts, in which the standardized Bernoulli distribution $\bs(p)=p\,\de_{\sqrt{q/p}}+q\,\de_{-\sqrt{p/q}}$ is replaced by the standardized symmetric three-point distribution 
$\ST(p):=\frac p2\,\de_{\,\sqrt{1/p}}
+(1-p)\,\de_0+\frac p2\,\de_{-\sqrt{1/p}}$, where $\de_a$ is the Dirac probability measure concentrated at one point $a$ and $p\in(0,1]$.
\big(At that,
the class $\F3$ may be replaced by the larger, ``two-tail'' version $\FF3$.
Note also that, since the distribution $\bs(p)$ is symmetric, one can assume without loss of generality that all generalized moment functions $f$ involved are even; indeed, for a symmetric r.v. $X$, one has $\E f(X)=\E g(X)$, where $g(x):=\frac12\,(f(x)+f(-x))$.\big)
Then the exact lower bound $m_*(p)$ on $m$ gets replaced by another exact lower bound, whose general expression is however more complicated than expression \eqref{eq:m_*} for $m_*(p)$. 
In fact, the exact lower bound on $m$ (denoted here by $m_{\st,\,\mathsf{Schur}}(p)$) for the symmetric-case analogues of the Schur-concavity statements in Theorem~\ref{th:main} turns out to be strictly greater for some values of $p$ than the exact lower bound on $m$ (denoted here by $m_\st(p)$) for the symmetric-case analogues of the above ``asymmetric'' results other than  Schur-concavity.

A simple upper bound on $m_{\st,\,\mathsf{Schur}}(p)$ and hence on $m_\st(p)$ for all $p\in(0,\frac12]$ is given by inequality
\begin{equation}\label{eq:mtst}
m_\st(p)\le m_{\st,\,\mathsf{Schur}}(p)\le m_{\st,\,\mathsf{high}}(p):=m_*(r)=
\frac{5-3\sqrt{1-2p}-2p} {4(\sqrt{\frac p2}+1-\sqrt{1-2p}-p)}
\end{equation}
where the function $m_*$ is defined by \eqref{eq:m_*} and
$r:=(1-\sqrt{1-2p})/2$ is the root in $(0,\frac12]$ of equation 
$2r(1-r)=p$. 
This follows because the convolution of the standardized Bernoulli distributions $\bs(r)$ and $\bs(1-r)$ is 
a symmetric three-point distribution, which ascribes probability $1-p$ to $\{0\}$ and therefore can be obtained from the standardized symmetric three-point distribution $\ST(p)$ by simple re-scaling. 
One can see that $m_{\st,\,\mathsf{high}}(p)\sim\frac1{\sqrt{2p}}$ as $p\downarrow0$ and $m_{\st,\,\mathsf{high}}(\frac12)=1$; let also $m_{\st,\,\mathsf{high}}(p):=1$ for all $p\in(\frac12,1)$.

A remarkable and not so difficult to prove fact is that
\begin{equation}\label{eq:p>1/2} 
m_\st(p)=1\quad\forall p\in[p_*,1],\quad\text{where}\quad p_*:=\sqrt2-1=0.4142\dots. 
\end{equation}
More generally, it appears that  
\begin{align*} 
m_\st(p) 
	&=m_1(p)\ii{0<p<p_{0,1}}
	+m_0(p)\ii{p_{0,1}\le p<p_*}
	+\ii{p_*\le p<1},\quad\text{} \notag 
\end{align*}
where $p_{0,1}=0.3878\dots$, $m_1(p):=\frac12\,\sqrt{\frac{2-p}p}$, and $m_0(p):=\frac1{2\log_2 z}$, where $z$ is (for each $p\in(0,p_*)$)
the only root in the interval $(0,\sqrt2)$ of polynomial 
\begin{multline*}
(9 p^2-24 p+16) z^6+(-36 p^2+120 p-96) z^5+(36 p^2-216
   p+240) z^4\\
   +(-20 p^3+60 p^2+120 p-320) z^3+(72 p^3-168
   p^2+96 p+240) z^2\\
   +(-96 p^3+144 p^2-144 p-96) z+4 p^4+40 p^3-44
   p^2+48 p+16.
\end{multline*}
(Here, as usual, $\ii{\mathcal{A}}$ denotes the indicator of an assertion $\mathcal{A}$.)

On the other hand, one can show that a lower bound on $m_\st(p)$ and hence on $m_{\st,\,\mathsf{Schur}}(p)$ is given by the inequality
\begin{equation}\label{eq:m-lower} 
\begin{aligned} 
& m_{\st,\,\mathsf{Schur}}(p)\ge m_\st(p)\ge m_{\st,\mathsf{low}}(p):=\max(1,m_1(p),m_{\mathsf{low}}(p)) \\
	&=m_1(p)\ii{0<p<p_{{\mathsf{low}},1}}
	+m_{\mathsf{low}}(p)\ii{p_{{\mathsf{low}},1}\le p<p_*}
	+\ii{p_*\le p<1},\quad\text{} 
\end{aligned}
\end{equation}
where $p_{{\mathsf{low}},1}=0.3889\dots$ and 
$m_{\mathsf{low}}(p):=\frac3{2(1+\log_2(1+p))}$. 
It appears that this lower bound, $m_{\st,\mathsf{low}}(p)$ is very close to $m_\st(p)$ and differs from the latter only for $p$ in the rather narrow interval $(p_{0,1},p_*)\approx
(0.3878,0.4142)$, and at most by $m_\st(p_{{\mathsf{low}},1})-m_{\st,\mathsf{low}}(p_{{\mathsf{low}},1})=0.0008598\dots$. 
It follows from \eqref{eq:m-lower} that the lower bound $p_*=\sqrt2-1$ on $p$ in \eqref{eq:p>1/2} is exact; that is, $m_\st(p)>1$ $\forall p\in(0,p_*)$.
Also, the four bounds, $m_{\st,\mathsf{low}}(p)\le m_\st(p)\le m_{\st,\,\mathsf{Schur}}(p)\le m_{\st,\mathsf{high}}(p)$, all equal $1$ for all $p\in[\frac12,1)$ and are asymptotic to $1/\sqrt{2p}$ as $p\downarrow0$;
moreover, the difference between any two of these four bounds goes to $0$ as $p\downarrow0$, since it is easy to see that $m_{\st,\mathsf{high}}(p)-m_{\st,\mathsf{low}}(p)\to0$ as $p\downarrow0$.

In this ``symmetric" setting, the place of supermartingales gets taken by conditionally symmetric martingales (that is, martingales $(S_i)$ with $S_0=0$ a.s.\ and conditionally symmetric differences $X_i$), and ``bounded-asymmetry'' condition \eqref{eq:bent-asymm} in Corollary~\ref{cor:bent} gets replaced by ``bounded-kurtosis'' condition
\begin{equation}\label{eq:bent-symm}
\dfrac{b_i^2}{c_i^2}\le\dfrac1p\quad\forall i=1,\dots,n.
\end{equation}
(Note that the kurtosis of a standardized symmetric three-point r.v. $\SST\sim\ST(p)$ is $\frac1p$.)
Thus, one now has another -- ``symmetric" -- version of inequality \eqref{eq:bent-Th1.3}, which immediately implies a striking ``symmetric'' version of inequality \eqref{eq:bent-Th1.3,normal}; namely, in view of \eqref{eq:p>1/2} and \eqref{eq:bent-symm}, inequality \eqref{eq:bent-Th1.3,normal} holds for all conditionally symmetric martingales if the condition $b_i=c_i$ $\forall i$ of \cite[Theorem~1.3]{bent-ap} is replaced by the much less restrictive condition
$$b_i^2\le\frac{c_i^2}{\sqrt2-1}\approx2.4142\,c_i^2\quad\forall i=1,\dots,n.$$

In particular, one has the following generalizations of \eqref{eq:khin} and \eqref{eq:pin94}, respectively:
\begin{equation}\label{eq:khin-st}
\E f\lp a_1\SST_1+\dots+a_n\SST_n\rp\le\E f(Z)\quad\forall f\in\FF3
\end{equation}
and
\begin{equation}\label{eq:pin94-st}
\P\lp a_1\SST_1+\dots+a_n\SST_n\ge x\rp\le\frac{2e^3}9\,\P(Z\ge x)\quad
\forall x\in\R, 
\end{equation}
where $\SST_i\overset{\text{i.i.d.}}\sim\ST(p)$ and $p\in[\sqrt2-1,1]$; 
note that $\ST(p)$ coincides with the distribution of a Rademacher r.v. if $p=1$.

Inequalities for the much simpler case when $b_i=b$ $\forall i$ were obtained in \cite{berger,bent-symm-exp} (in the ``exponential'' case) and in \cite{bent-symm}. 

One also has multi-dimensional analogues of inequalities \eqref{eq:khin-st} and \eqref{eq:pin94-st}. Namely, for all $p\in[\sqrt2-1,1]$ one has the following generalizations of inequalities \eqref{eq:khin-hilbert} and \eqref{eq:pin94-hilbert}:
\begin{equation}\label{eq:khin-st-hilbert}
\E f\lp \|\SST_1\xx_1+\dots+\SST_n\xx_n\| \rp\le\E f(|Z|)\quad
\text{for all even $f$ in $\FF3$}
\end{equation}
and
\begin{equation}\label{eq:pin94-st-hilbert}
\P\lp  \|\SST_1\xx_1+\dots+\SST_n\xx_n\| \ge x \rp\le\frac{2e^3}9\,\P(|Z|\ge x)\quad
\forall x\in\R, 
\end{equation}
with the Rademacher r.v.'s $\vp_1,\dots,\vp_n$ replaced by $\SST_1,\dots,\SST_n$.

Condition~\eqref{eq:bent-symm} may be interpreted as a condition of boundedness of the kurtoses of the $X_i$'s; cf. bounded-asymmetry conditions \eqref{eq:bi/ai} and \eqref{eq:Bi/Ai}. 
Note that the usual interpretation of the term ``kurtosis" (as well as the term itself, meaning ``peakedness") is not quite adequate, especially in contexts of large deviations. Indeed, without loss of generality, let us assume that a r.v. $X$ is standardized, that is, zero-mean and unit-variance. 
Then the kurtosis of $X$ is 
$\kurt(X)=\E X^4=1+\E(X^2-1)^2=1+\Var(X^2)$,
which clearly is the measure of relative spread of the values of the nonnegative r.v. $X^2$ about its mean $1$. Thus, $\kurt(X)$ is large if and and only if at least one of the two tails, left or right, of the distribution of $X^2$ on $[0,\infty)$ is heavy. These two tails can be measured, respectively, by the ``negative" and ``positive'' parts of $\kurt(X)-1$, namely $\kurt_-(X):=\E(X^2-1)^2\ii{X^2<1}$ and $\kurt_+(X):=\E(X^2-1)^2\ii{X^2>1}$. 
If $\kurt_-(X)$ is large, then the distribution of $X$ is peaked at $0$; 
if $\kurt_+(X)$ is large, then at least one tail of the distribution of $X$ is heavy. Of these two reasons for the kurtosis of $X$ to be large, the heaviness of the tails of the distribution of $X$ seems more important than the peakedness, at least where large deviations are involved. Also, note that the kurtosis of $X$ can be infinite only if the tails of the distribution of $X$ are heavy enough. As was noted, the kurtosis of a standardized symmetric three-point r.v. $\SST\sim\ST(p)$ is $\frac1p$, which is decreasing in $p$; at that, both the negative part $\kurt_-(\SST)=1-p$ and the positive part $\kurt_+(\SST)=4\sinh^2\ln\frac1{\sqrt p}$ are decreasing in $p$.
More generally, it is not difficult to show that for any functions $f\in\FF3$ (except for $f$ of the form $f(x)\equiv a+bx^2$ for some nonnegative real constants $a$ and $b$) the generalized moments $\E f(\SST)$ and hence $\E f(c_1\SST_1+\dots+c_n\SST_n)$ are strictly decreasing in $p\in(0,1]$, where $\SST,\SST_i\overset{\text{i.i.d.}}\sim\ST(p)$ and the $c_i$'s are nonzero real constants. 
Thus, inequality \eqref{eq:khin-st} 
for $p\in[\sqrt2-1,1)$
not only generalizes but also strengthens inequality \eqref{eq:khin}.
Similarly, inequality \eqref{eq:pin94-st} is stronger than inequality 
\eqref{eq:pin94} for $p\in[\sqrt2-1,1)$ and large enough $x$.

It may also be instructive to compare inequality \eqref{eq:pin94-st} with available asymptotic results. For instance, it follows from \cite{feller} or \cite{petrov} that for \emph{every} $p\in(0,1)$
\begin{equation}\label{eq:asymp-st}
\P\lp \SST_1+\dots+\SST_n\ge x\rp\sim\P(Z\ge x/\sqrt n)
\end{equation}
provided that $n\to\infty$ and $x=o(n)$. The advantages of this asymptotics over inequality \eqref{eq:pin94-st} are (i) that \eqref{eq:asymp-st} holds for every $p\in(0,1)$ and not just for $p\in[\sqrt2-1,1]$ and (ii) there is no ``extra'' constant factor (such as $2e^3/9$) in the right-hand side of \eqref{eq:asymp-st}. However, these advantages seem to be counterbalanced by the following: (i) the asymptotic relation \eqref{eq:asymp-st}, without an upper bound on the rate of convergence, is impossible to use in statistical practice when one needs to be certain that the left-hand side of \eqref{eq:asymp-st} does not exceed a prescribed level; (ii) the asymptotics is valid only in the zone $x=o(n)$, and this zone is defined only qualitatively; (iii) the summands $\SST_1,\dots,\SST_n$ in the left-hand side of \eqref{eq:asymp-st} are identically distributed; if coefficients $a_1,\dots,a_n$ are allowed in \eqref{eq:asymp-st} (as in \eqref{eq:pin94-st}), they cannot differ too much from another, and also then the zone $x=o(n)$ must be narrowed; (iv) there is no ``extra'' constant factor (such as $2e^3/9$) in the right-hand side of the generalized-moment comparison inequality \eqref{eq:khin-st}; (v) one can use the more precise upper bounds of the form given in \eqref{eq:comp-prob2} or \eqref{eq:interp} instead of \eqref{eq:pin94-st}. 

Details on this long remark, Remark~\ref{rem:symm}, will be presented elsewhere.
\end{remark}

\begin{remark}\label{rem:symm-exp}
Directions of further research indicated in Remarks~\ref{rem:exp} and \ref{rem:symm} can be combined. That is, one can at once replace (i) the class of generalized moment functions $\F3$ by the much narrower class $\F\infty$ (the latter being in effect the same as the class of all nondecreasing exponential functions) and (ii) the standardized Bernoulli distribution by the standardized symmetric three-point distribution. At that, one has another exact lower bound, say $m_{\mathsf{st,\,exp}}(p)$, on $m$. 
Obviously, $m_{\mathsf{st,\,exp}}(p)\le m_\st(p)$ for all $p\in(0,1]$. 
In particular, $m_{\mathsf{st,\,exp}}(p)=m_\st(p)=1$ for all $p\in[\sqrt{2}-1,1]$. 
Also, similarly to \eqref{eq:mtst} and in view of \eqref{eq:m-exp-up}, one has 
$m_{\mathsf{st,\,exp}}(p)\le m_{\mathsf{exp}}(r)\le m_{\mathsf{exp,\,up}}(r)$ for all $p\in(0,\sqrt{2}-1)$ (and even for all $p\in(0,\frac12]$),
where as before $r=(1-\sqrt{1-2p})/2$ is the root in $(0,\frac12]$ of equation $2r(1-r)=p$. 
While, as noted, $m_{\mathsf{st,\,exp}}(p)=1$ for all $p\in[\sqrt{2}-1,1]$, one has $m_{\mathsf{st,\,exp}}(p)>1$ for all $p\in(0,\frac13)$. 
Indeed, if $m_{\mathsf{st,\,exp}}(p)=1$ for some $p\in(0,1)$, then 
inequality \eqref{eq:khin-st} holds for that same $p$ and all $f\in\F\infty$. Taking now $n=1$ and $f(x)\equiv e^{kx}$ in \eqref{eq:khin-st}, and then letting $k\downarrow0$, one has $\frac1p=\E\SST_1^4\le\E Z^4=3$, whence $p\ge\frac13$.

Details on this remark, Remark~\ref{rem:symm-exp}, will be presented elsewhere.
\end{remark}

\subsection{Bounds on self-normalized sums}\label{t-stats} 
(Details on the results presented in this subsection will be given elsewhere.)
Efron \cite{efron} considered the so-called self-normalized sum
\begin{equation}\label{eq:V}
V:=\frac{X_1+\dots+X_n}{\sqrt{X_1^2+\dots+X_n^2}},
\end{equation}
assuming that the $X_i$'s satisfy the orthant symmetry condition:
the joint distribution of $\de_1X_1,\dots,\de_n X_n$ is the same for any choice of signs $\de_1,\dots,\de_n\in\{1,-1\}$, so that, in particular, each $X_i$ is symmetric(ally distributed). It suffices that the $X_i$'s be independent and symmetrically (but not necessarily identically) distributed. 
On the event $\{X_1=\dots=X_n=0\}$, let $V:=0$.

In Subsection~\ref{t-stats} we assume throughout that the $X_i$'s are all non-degenerate: $\P(X_i=0)<1$ for all $i$.

Note that the conditional distribution of any symmetric r.v.\ $X$ given $|X|$ is the symmetric distribution on the (at most) two-point set $\{X,-X\}$. 
Therefore, under the orthant symmetry condition, the distribution of $V$
is the mixture of the distributions of the normalized Khinchin-Rademacher sums 
$a_1\vp_1+\dots+a_n\vp_n$, where $a_i=X_i/(X_1^2+\dots+X_n^2)^{\frac12}$, so that $a_1^2+\dots+a_n^2=1$ (except on the event $\{X_1=\dots=X_n=0\}$, where $a_1=\dots=a_n=0$).
Hence, by \eqref{eq:khin} (for $f(x)\equiv e^{\la x}$, $\la\ge0$) and \eqref{eq:exp} one has 
\begin{equation}\label{eq:khin-V}
\E e^{\la V}\le\E e^{\la Z}\quad\forall\la\ge0
\end{equation}
and
\begin{equation}\label{eq:exp-V}
\P\lp V\ge x\rp\le e^{-x^2/2}\quad \forall x\ge0.
\end{equation}
These results can be easily restated in terms of Student's statistic $T$, which is a monotonic function of $V$, as noted by Efron;
namely, $T=\sqrt{\frac{n-1}n}\,V/\sqrt{1-V^2/n}$.
Inequalities \eqref{eq:khin-V} and \eqref{eq:exp-V} were improved in \cite{eaton1,eaton2,pin94} using \eqref{eq:khin} (for $f\in\H3$) and \eqref{eq:pin94}, so that one has 
\begin{equation}\label{eq:khin-V-pin94}
\E f(V)\le\E f(Z)\quad\forall f\in\H3
\end{equation}
and
\begin{equation}\label{eq:exp-V-pin94}
\P\lp V\ge x\rp\le\frac{2e^3}9\,\P(Z\ge x)\quad
\forall x\in\R.
\end{equation}
Multivariate analogues of these results, which can be expressed in terms of Hotelling's statistic in place of Student's, were also obtained in \cite{pin94}. 

It was pointed out in \cite[Theorem~2.8]{pin94} that, since the normal tail decreases fast, inequality \eqref{eq:exp-V-pin94} implies that relevant quantiles of $V$ may exceed the corresponding standard normal quantiles only by a relatively small amount, so that one can use \eqref{eq:exp-V-pin94} rather efficiently to test symmetry even for non-i.i.d.\ observations. 

Here we shall present extensions of inequalities \eqref{eq:khin-V-pin94} and \eqref{eq:exp-V-pin94} to the case when the $X_i$'s are not symmetric, as well as improvements of \eqref{eq:khin-V-pin94} and \eqref{eq:exp-V-pin94} in the ``symmetric case''. Asymptotics for large deviations of $V$ for i.i.d.\ $X_i$'s without moment conditions was obtained recently by Jing, Shao and Zhou \cite{jing-shao}. Comments similar to the ones made at the end of Remark~\ref{rem:symm} apply here as well.

\subsubsection{The ``asymmetric'' case} \label{V-asymm}
The basic idea here is to represent any zero-mean, possibly asymmetric distribution as an appropriate mixture of two-point zero-mean distributions.  
Let us assume at first that a zero-mean r.v.\ $X$ has an everywhere strictly positive density function. Consider the truncated r.v.\ $X_{a,b}:=X\ii{a\le X\le b}$. 
Then, for every fixed $a\in(-\infty,0]$, the function $b\mapsto\E X_{a,b}$ is continuous and increasing on the interval $[0,\infty)$ from $\E X_{a,0}\le0$ to 
$\E X_{a,\infty}>0$. Hence, for each $a\in(-\infty,0]$, there exists a unique value $b\in[0,\infty)$ such that $\E X_{a,b}=0$. Similarly, for each $b\in[0,\infty)$, there exists a unique value $a\in(-\infty,0]$ such that 
$\E X_{a,b}=0$. That is, one has a one-to-one correspondence between $a\in(-\infty,0]$ and $b\in[0,\infty)$ such that $\E X_{a,b}=0$. Denote by 
$\r:=\r_X$ the {\em reciprocating} function defined on $\R$ and carrying this correspondence, so that 
$$\E X\ii{\text{$X$ is between $x$ and $\r(x)$} }=0\quad\forall x\in\R,$$
the function $\r$ is decreasing on $\R$ and such that $\r(\r(x))=x$ $\forall x\in\R$; moreover, $\r(0)=0$. 
(Clearly, $\r(x)=-x$ for all real $x$ if the r.v.\ $X$ is symmetric.)
One also has
\begin{equation} \label{eq:hat}
\r(x)=
x_-(G(x))\,\ii{x>0}+
x_+(G(x))\,\ii{x<0}, 
\end{equation}
where $x_{\pm}(h)$ stand for the positive and negative roots $x$ of the equation $G(x)=h$ and, in turn,
\begin{equation}\label{eq:G}
G(x):=
\E|X|\ii{|X|\le|x|,\;\sign X=\sign x}.
\end{equation}
Thus, the set 
$\{\,\{x,\r(x)\}\colon x\in\R\,\}$ of (at-most-)two-point sets constitutes a partition of $\R$. Moreover, the two-point set $\{x,\r(x)\}$ is uniquely determined by the distance $|x-\r(x)|$ between the two points, as well as by the product $|x|\,|\r(x)|$. 
Now one can see that the conditional distribution of the zero-mean r.v.\ $X$ given $W:=|X-\r(X)|$ (or, equivalently, $Y:=|X\,\r(X)|$) is the uniquely determined zero-mean distribution on the two-point set $\{X,\r(X)\}$.
Thus, the distribution of the zero-mean r.v.\ $X$ with an everywhere positive density is represented as a mixture of two-point zero-mean distributions.
This mixture is given rather explicitly, provided that the distribution of r.v.\ $X$ is known. 

Thus, one has generalized versions of the self-normalized sum \eqref{eq:V}, which require -- instead of the symmetry of independent r.v.'s $X_i$ -- only that the $X_i$'s be zero-mean:
\begin{equation} \label{eq:V_W} 
V_W:=\frac{X_1+\dots+X_n}{\frac12\sqrt{W_1^2+\dots+W_n^2}}\quad\text{and}\quad
V_{Y,m}:=\frac{X_1+\dots+X_n}{(Y_1^m+\dots+Y_n^m)^{\frac1{2m}}},
\end{equation}
where $m\ge1$,
$$W_i:=|X_i-\r_i(X_i)|,\quad\text{and}\quad Y_i:=|X_i\,\r_i(X_i)|,$$
and the reciprocating function $\r_i:=\r_{X_i}$ is constructed as above, based on the distribution of $X_i$, for each $i$, so that the reciprocating functions $\r_i$ may be different from one another if the $X_i$'s are not identically distributed.
On the event $\{X_1=\dots=X_n=0\}$ (which is the same as either one of events $\{W_1=\dots=W_n=0\}$ and $\{Y_1=\dots=Y_n=0\}$), let $V_W:=0$ and $V_{Y,m}:=0$. 
Note that $V_W=V_{Y,1}=V$ when the $X_i$'s are symmetric.
Logan \emph{et al} \cite{logan} and Shao \cite{shao} obtained limit theorems for the ``symmetric'' version of $V_{Y,m}$ (with the reciprocating function $\r(x)\equiv-x$), whereas the $X_i$'s did not need to be symmetric.

These constructions can be extended to the general case of any zero-mean r.v.\ $X$, absolutely continuous or not. Here, one can use randomization (by means of a r.v.\ uniformly distributed in interval $(0,1)$) to deal with the atoms of the distribution of r.v.\ $X$, and a modification of the inverse functions $x_{\pm}(h)$ to deal with the intervals on which the distribution function of $X$ and hence the function $G$ are constant. Namely, in general $\r(X)$ is replaced by $\r(X,U)$, where $U$ is a r.v.\ uniformly distributed in interval $(0,1)$ and independent of $X$ and,
for $x\in\R$ and $u\in(0,1)$,  
\begin{gather*} 
\r(x,u):=
\begin{cases}
x_-(G(x-)+u\cdot(G(x)-G(x-))) & \text{ if }x\in[0,\infty), \\
x_+(G(x+)+u\cdot(G(x)-G(x+))) & \text{ if }x\in(-\infty,0],
\end{cases}
\\
\begin{aligned}
x_+(h) & := \inf\{x\in[0,\infty]\colon G(x)\ge h\},  \\
x_-(h)& :=\sup\{x\in[-\infty,0]\colon G(x)\ge h\}. 
\end{aligned}
\end{gather*}

By conditioning on the $W_i$'s or $Y_i$'s one obtains the following corollaries.

\begin{corollary} \label{cor:stud-normal}
\emph{(From results of \cite{normal}:)} 
\begin{align} 
\E f(V_W) &\le \E f(Z)\quad\forall f\in\H5,\quad\text{whence} \label{eq:stud-f(Z)} \\
\P(V_W\ge x)&\le c_{5,0}\P(Z\ge x)\quad\forall x\in\R, \label{eq:stud-P(Z>x)}
\end{align}
where, in accordance with \eqref{eq:c(al,beta)}, $c_{5,0}=5!(e/5)^5=5.699\dots$. 
\end{corollary}  

\begin{corollary} \label{cor:stud-asymm}
\emph{(From Theorem~\ref{th:bern} and Corollary~\ref{cor:beta,prob}:)} 
Suppose that for some $p\in(0,1)$ and all $i\in\{1,\dots,n\}$
\begin{equation}\label{eq:bounded-asymm}
\frac{X_i}{|\r(X_i)|}\,\ii{X_i>0} \le \frac qp\ \text{a.s.}	
\end{equation}
Then for all $m\ge m_*(p)$
\begin{align} 
\E f(V_{Y,m}) &\le \E f(T_n)\quad\forall f\in\F3\quad\text{and} \label{eq:stud-f(T_n)} \\
\P(V_{Y,m}\ge x)&\le c_{3,0}\P^\lc(T_n\ge x)\quad\forall x\in\R, \label{eq:stud-P(T_n>x)}
\end{align}
where $T_n$ and $\P^\lc(T_n\ge x)$ have the same meaning as in Corollary~\ref{cor:beta,prob} and, in accordance with \eqref{eq:c(al,beta)}, $c_{3,0}=2e^3/9=4.4634\ldots$. 
\end{corollary}  

Of course, one can replace the upper bound in inequalities like \eqref{eq:stud-P(T_n>x)} by either of the more precise (but slightly less transparent and more difficult to compute) upper bounds given in \eqref{eq:prob2} and \eqref{eq:prob5}. 

By Remark~\ref{rem:exp}, inequality \eqref{eq:stud-f(T_n)} will continue to hold for all $m\ge m_{\mathsf{exp}}(p)$ provided that the class $\F3$ is replaced by the smalled class $\F\infty$; correspondingly, for such $m$ one will have the exponential upper bound on $\P(V_{Y,m}\ge x)$ of the form $e^{-nH}$ as in \eqref{eq:prob4}. 

Condition \eqref{eq:bounded-asymm} is likely to hold when the $X_i$'s are bounded i.i.d.\ r.v.'s. 

Note that the reciprocating function $\r$ depends on the (usually unknown in statistics) distribution of the underlying r.v.\ $X$. However, if e.g.\ the $X_i$ constitute an i.i.d.\ sample, then the function $G$ defined by \eqref{eq:G} can be estimated based on the sample, so that one can estimate the reciprocating function $\r$. Thus, replacing $X_1+\dots+X_n$ in the numerators of $V_W$ and $V_{Y,m}$ by $X_1+\dots+X_n-n\theta$,
one obtains approximate pivots to be used to construct confidence intervals or, equivalently, tests for an unknown mean $\theta$. 
One can also use bootstrap to estimate the distributions of such pivots.

\subsubsection{The ``symmetric'' case} \label{V-symm}
Here we assume that $X_1,\dots,X_n$ are independent symmetric r.v.'s. In this case, we already have the upper bounds given by \eqref{eq:khin-V-pin94} and \eqref{eq:exp-V-pin94}, which are based on \eqref{eq:khin} and \eqref{eq:pin94}.
As was pointed out, bounds \eqref{eq:khin-st} and \eqref{eq:pin94-st} improve and generalize \eqref{eq:khin} and \eqref{eq:exp}. Correspondingly, the upper bounds given by \eqref{eq:khin-V-pin94} and \eqref{eq:exp-V-pin94} can be improved and generalized as follows. 

Introduce another, ``symmetric'' modification of the standard self-normalized sum $V$ given by \eqref{eq:V}:
\begin{equation}\label{eq:V-symm}
V_{\mathsf{symm},m,p}:=\frac{\SST_1X_1+\dots+\SST_n X_n}
{(|X_1|^{2m}+\dots+|X_n|^{2m})^{\frac1{2m}}},
\end{equation}	
where $p\in(0,1)$, $m\ge1$, $\SST_1,\dots,\SST_n$ are i.i.d.\ $\ST(p)$ r.v.'s, independent also of $X_1,\dots,X_n$.
On the event $\{X_1=\dots=X_n=0\}$, let $V_{\mathsf{symm},m,p}:=0$.

Thus, the distribution of $V_{\mathsf{symm},m,p}$
is the mixture of the distributions of the linear combinations of the form 
$a_1\SST_1+\dots+a_n\SST_n$, where $a_i=X_i/(|X_1|^{2m}+\dots+|X_n|^{2m})^{\frac1{2m}}$, so that $|a_1|^{2m}+\dots+|a_n|^{2m}=1$ (except on the event $\{X_1=\dots=X_n=0\}$, where $a_1=\dots=a_n=0$). 

By Remark~\ref{rem:symm}, it follows that
for all $m\ge m_\st(p)$
\begin{align} 
\E f(V_{\mathsf{symm},m,p}) &\le \E f(n^{-1/(2m)}(\SST_1+\dots+\SST_n))\quad\forall f\in\FF3,\quad\text{whence} \label{eq:stud-f(V-symm-1)} \\
\P(V_{\mathsf{symm},m,p}\ge x)&\le c_{3,0}\P^\lc(n^{-1/(2m)}(\SST_1+\dots+\SST_n)\ge x)\quad\forall x\in\R, \label{eq:stud-P(V-symm-1>x)}
\end{align}
where $\P^\lc$ again denotes the least log-concave majorant of the corresponding tail function. 
In particular (cf. \eqref{eq:khin-st} and \eqref{eq:pin94-st}), for all $p\in[\sqrt2-1,1)$,  
\begin{align} 
\E f(V_{\mathsf{symm},1,p}) 
&\le \E f(n^{-1/2}(\SST_1+\dots+\SST_n))\\
&\le \E f(Z)\quad\forall f\in\FF3, \label{eq:stud-f(V-symm-1-Z)} \\
\P(V_{\mathsf{symm},1,p}\ge x)
&\le c_{3,0}\P^\lc(n^{-1/2}(\SST_1+\dots+\SST_n)\ge x) \\
&\le c_{3,0}\P(Z\ge x)\quad\forall x\in\R. \label{eq:stud-P(V-symm-1>x)-Z}
\end{align}

In view of the mixture representations of the distributions of $V$ and $V_{\mathsf{symm},m,p}$ and the discussion after inequalities \eqref{eq:khin-st}--\eqref{eq:pin94-st-hilbert}, inequalities \eqref{eq:stud-f(V-symm-1)}--\eqref{eq:stud-P(V-symm-1>x)-Z} for $V_{\mathsf{symm},m,p}$ generalize/improve inequalities \eqref{eq:khin-V}--\eqref{eq:exp-V-pin94} for $V$.

The classic self-normalized sum $V$ can be obviously used, e.g., as a test statistic to test the symmetry of the distributions of the $X_i$'s.
However, it is not seen how $V_{\mathsf{symm},m,p}$ can be used in symmetry tests, because its distribution will be always symmetric, even if the distributions of the $X_i$'s are not. 

On the other hand, there are two issues with $V$: 
\begin{description}
	\item[(i)] 
	$V$ may have too light tails if the tails of the $X_i$'s are heavy enough, and so, the symmetry test based on such an upper bound as the one given by inequality \eqref{eq:exp-V-pin94} may turn out to be too conservative and hence lacking some power;
	\item[(ii)] 
	the way $V$ (as well as its modifications considered so far -- $V$, $V_W$, $V_{Y,m}$, and $V_{\mathsf{symm},m,p}$) ``deals'' with the event when all the $X_i$'s take on the zero value certainly seems to be too conservative; note that this event may naturally occur with a nonzero probability if the distribution of $X_i$ is discrete or if the original r.v.'s $X_i$ are replaced by the corresponding truncated r.v.'s $X_i\ii{|X_i|<b_i}$ for some $b_i>0$ (which may be done to 
increase the power of the test). 
However, note that usually the probabilities $p_i:=\P(X_i\ne0)$ will be close to (even if less than) $1$. 
\end{description}

To try to resolve these two issues with $V$, we shall suggest yet another modification of it.
Each of the other self-normalized sums introduced above -- $V_W$, $V_{Y,m}$, and $V_{\mathsf{symm},m,p}$ -- can be modified in the same manner.  

Let us accompany any r.v.\ $X$ with a r.v.\ of the form
\begin{equation}\label{eq:hat X}
\hat X:=X+\tilde X\,\ii{X=0},
\end{equation}
where $\tilde X$ is a r.v.\ which is independent of $X$ and whose distribution coincides with the conditional distribution of $X$ given that $X\ne0$; thus, $\hat X=X$ if $X\ne0$ and $\hat X=\tilde X$ if $X=0$, so that $\P(\hat X\ne0)=1$. It is not hard to see that, if $X$ is symmetric, then the conditional distribution of $X$ given $|\hat X|$ coincides with the symmetric distribution 
$\frac p2\de_{\hat X}+(1-p)\de_0+\frac p2\de_{-\hat X}$ on the three-point set $\{\hat X,0,-\hat X\}$, where $p:=\P(X\ne0)$:
\begin{equation}\label{eq:cond}
\mathcal{L}(X|\hat X)=\tfrac p2\de_{\hat X}+(1-p)\de_0+\tfrac p2\de_{-\hat X}.
\end{equation}
Moreover, the distribution of $\hat X$ coincides (just as that of $\tilde X$ does) with the conditional distribution of $X$ given that $X\ne0$.
Indeed, $\hat X$ is equal in distribution to
the first nonzero member of an infinite random sequence $(X^{(1)},X^{(2)},\dots)$, 
if $X,X^{(1)},X^{(2)},\dots$ are i.i.d.\ r.v.'s; therefore, one may assume that 
$\hat X$ equals (as a r.v., and not just in distribution) to
the first nonzero member of the sequence $(X,X^{(1)},X^{(2)},\dots)$.
Thus, roughly speaking, to get $\hat X$, one samples from the distribution of $X$ till getting a nonzero replica of $X$.
Note also that for any even function $g$ such that $g(0)=0$ one has
\begin{equation}\label{eq:var}
\Var g(\hat X)-\Var g(X)
=\tfrac{1-p}p\,\big(\Var g(X)-\tfrac1p(\E g(X))^2\big).
\end{equation}

It follows from \eqref{eq:cond} that the distribution of yet another modification of the self-normalized sum $V$,
\begin{equation}\label{eq:V-symm-2}
\hat V_{\mathsf{symm},m,p}:=\frac{X_1+\dots+X_n}
{\sqrt{p}\,(|\hat X_1|^{2m}+\dots+|\hat X_n|^{2m})^{\frac1{2m}}},
\end{equation}
is the mixture of the distributions of the linear combinations of the form \\
$a_1\SST_1(p_1)+\dots+a_n\SST_n(p_n)$, where the $\SST_i(p_i)$'s are independent r.v.'s such that $\SST_i(p_i)\sim\ST(p_i)$ for all $i$ and
$a_i=\hat X_i/(|\hat X_1|^{2m}+\dots+|\hat X_n|^{2m})^{\frac1{2m}}$, so that $|a_1|^{2m}+\dots+|a_n|^{2m}=1$ a.s.; here, $p\in(0,1)$,
$$p\ge p_i:=\P(X_i\ne0)\quad\forall i,$$ 
and each $\hat X_i$ is produced based on $X_i$ according to formula \eqref{eq:hat X}, where the $\tilde X_i$'s are independent of one another and of the $X_i$'s.

Recall that, for each $i$, one has $\P(\hat X_i=0)=0$ and the distribution of $\hat X_i$ coincides with the conditional distribution of $X_i$ given that $X_i\ne0$. 

It follows that 
\begin{align} 
\E f(\hat V_{\mathsf{symm},m,p}) &\le \E f(n^{-1/(2m)}(\SST_1+\dots+\SST_n))\quad\forall f\in\FF3, \label{eq:stud-f(V-symm-2)} \\
\intertext{whence}
\P(\hat V_{\mathsf{symm},m,p}\ge x)&\le c_{3,0}\P^\lc(n^{-1/(2m)}(\SST_1+\dots+\SST_n)\ge x)\quad\forall x\in\R, \label{eq:stud-P(V-symm-2>x)}
\end{align}
provided that
$$
m\ge m_\st(p);$$
in particular, if $p\ge\sqrt2-1$ (which will typically be the case),  
\begin{align} 
\E f(\hat V_{\mathsf{symm},1,p}) 
&\le \E f(n^{-1/2}(\SST_1+\dots+\SST_n))\\
&\le \E f(Z)\quad\forall f\in\FF3, \label{eq:stud-f(V-symm-2-Z)} \\
\P(\hat V_{\mathsf{symm},1,p}\ge x)
&\le c_{3,0}\P^\lc(n^{-1/2}(\SST_1+\dots+\SST_n)\ge x) \label{eq:stud-P(V-symm-2>x)-T_n} \\
&\le c_{3,0}\P(Z\ge x)\quad\forall x\in\R; \label{eq:stud-P(V-symm-2>x)-Z}
\end{align}
here again the $\SST_i$'s are i.i.d.\ $\ST(p)$ r.v.'s. 

Note that $\esssup\hat V_{\mathsf{symm},1,p}=\sqrt{n/p}$, while $\esssup V=\sqrt{n}$. Thus, the tails of $\hat V_{\mathsf{symm},1,p}$ are longer in some sense than those of $V$, so that the symmetry test based on an inequality such as \eqref{eq:stud-P(V-symm-2>x)-T_n} or \eqref{eq:stud-P(V-symm-2>x)-Z} may be less conservative and hence more powerful than the corresponding test based on an inequality such as \eqref{eq:exp-V-pin94}, especially if the size of the test is small enough. 
An interesting question is how the generalized moments $\E f(\hat V_{\mathsf{symm},m,p})$ compare with $\E f(V)$ for $f\in\FF3$. 

It may also be of interest to compare $\E f(\hat V_{\mathsf{symm},m,p})$
with
$\E f(V_{\mathsf{symm},m,p})$ (assuming that $p_i=p$ for all $i$). It seems that neither of them dominates the other one in general. In view of \eqref{eq:var}, it seems likely that $\E f(\hat V_{\mathsf{symm},m,p})$ will be greater than $\E f(V_{\mathsf{symm},m,p})$ if the tails of the distributions of the $X_i$'s are not too heavy: then the variability of the $|\hat X_i|$'s will be less than that of the $|X_i|$'s and hence, heuristically, the tails of $\hat V_{\mathsf{symm},m,p}$ will be heavier than those of $V_{\mathsf{symm},m,p}$.

Anyway, a definite advantage of $\hat V_{\mathsf{symm},m,p}$ (over $V$ and $V_{\mathsf{symm},m,p}$) is that the denominator of its ratio expression in \eqref{eq:V-symm-2} is nonzero a.s. 
On the other hand, an obvious disadvantage of $\hat V_{\mathsf{symm},m,p}$ is that to compute its value one needs to know the distributions of the $X_i$'s. 
In statistical practice, $\hat V_{\mathsf{symm},m,p}$ may be approximated (at least in the case when the $X_i$'s are i.i.d.) by replacing 
$\tilde X_i$ for each $i$ by $X_{J(i)}$, where $J(1),\dots,J(n)$ are (conditionally, given $X_1,\dots,X_n$) i.i.d.\ r.v.'s, each $J(i)$ having the uniform distribution on the set 
$$\{j\in\{1,\dots,n\}\colon X_j\ne0\};$$
this set is nonempty with a probability close to $1$ if $p$ is close to $1$ or $n$ is large; however, if this set happens to be empty, one can just set $\tilde X_i:=0$, so that one has $\hat X_i=0$. It would be interesting to compare generalized moments and tails of this ``practical'' version of $\hat V_{\mathsf{symm},m,p}$ with those of $\hat V_{\mathsf{symm},m,p}$ itself.

\section{Proofs}\label{proofs} 

\subsection{Statements of lemmas and proofs of the main results}\label{proofs1} 

The proofs of the main results are preceded in this subsection by some definitions and a series of lemmas. At least one of them (Lemma \ref{prop:superm}) may be of independent interest. The proofs of the lemmas are deferred further to Subsection \ref{proofs2}.

Let us introduce more classes of functions, in addition to the classes $\H3$, $\F3$, 
$\Fminus3$, 
$\FF3$, and $\GG3$ (recall \eqref{eq:H}, \eqref{eq:F3}, \eqref{eq:Fminus3}, \eqref{eq:FF3}, and \eqref{eq:GG3}:
\begin{align} 
\G3 & :=
\{f\colon\exists a\in\R,b\in\R,h\in\H3\ \forall x\in\R\ f(x)=a+b\,x+h(x)\}; 
\label{eq:G3}\\ 
\GPlus3 & :=
\{f\colon\exists a\in\R,b\ge0,h\in\H3\ \forall x\in\R\ f(x)=a+b\,x+h(x)\}. 
\label{eq:GPlus3}
\end{align} 

\begin{remark}\label{rem:shift}
It is not difficult to see that, if a function $f$ is in $\F3$ or any other defined above class of functions, then the shifted function 
$x\mapsto f(x+a)$ is also in the same class, for any real constant $a$. 
That is, all these classes of functions are shift-invariant. 
\end{remark}

\begin{lemma} \label{lem:g',g'',g'''}
Suppose that a function $g\colon\R\to\R$ is convex and such that there exists a finite limit $g(-\infty):=\lim_{x\to-\infty}g(x)$; in particular, the latter condition will obviously be the case if $g$ is nonnegative and nondecreasing.
Then 
$g'(-\infty)=0$, where $g'$ is the right derivative of $g$. 
\end{lemma}

\begin{lemma}\label{lem:F-al}
If $\al$ is a natural number then $\H\al$ coincides with the class 
$\tilde{\H\al}$ of all functions $f\colon\R\to\R$ such that
the derivative $f^{(\al-1)}$ is everywhere finite and convex, and $f^{(0)}(-\infty)=\dots=f^{(\al-1)}(-\infty)=0$. 
Moreover, if $f\in\H\al$, then all the functions $f^{(0)},\dots,f^{(\al-1)}$ are nonnegative. 
\end{lemma}

\begin{lemma} \label{lem:f-in-clG3}
Let $f\colon\R\to\R$ be a function such that $f''$ is finite, nonnegative, nondecreasing, and convex, with $f''(-\infty)=0$. Then $f\in\cl\G3$.
If, moreover, $f$ is nondecreasing, then $f\in\cl\GPlus3$.
\end{lemma}

\begin{lemma} \label{lem:F3=clG3}
One has $\F3=\cl\GPlus3$.
\end{lemma}

\begin{lemma} \label{lem:f3,ff3}
One has $\FF3=\cl\GG3$, where $\GG3$ is defined by \eqref{eq:GG3}.
\end{lemma}

\begin{lemma} \label{lem:f at infty}
If $f\in\F3$, then either $f(x)=O(x)$ as $x\to\infty$ or \\
$\liminf_{x\to\infty}f(x)/x^2\in(0,\infty]$. 
\end{lemma}

\begin{lemma} \label{lem:f at -infty}
If $f\in\F3$, then $f(x)=O(|x|)$ as $x\to-\infty$. 
\end{lemma}

\begin{lemma} \label{lem:FF3neGG3}
$\FF3\ne\GG3$.
\end{lemma}

\begin{proposition} \label{prop:F3neG3}
There exists a function $g\in\F3\setminus\G3$.
Since $\GPlus3\subseteq\G3$, it follows that
$\F3\ne\GPlus3$. {\em (This proposition complements Lemmas~\ref{lem:F3=clG3} and \ref{lem:FF3neGG3}; it will not be used elsewhere in this paper.) }
\end{proposition}

The following two lemmas are essentially well known. Their statements (and proofs) are given here for easy reference. 

\begin{lemma}\label{lem:extreme-bounded} \emph{(Cf. e.g. \cite{hunt} and \cite[Lemma~4.3]{bent-ap}.)}
Let $X$ be a r.v. such that $\E X\le0$ and $-a\le X\le b$ a.s. for some positive real numbers $a$ and $b$. 
Let $\BS\sim\bs(p)$ with
$$p:=\frac a{b+a}.$$ 
Then
$$\E f(X)\le\E f(\sqrt{a\,b}\,\BS)$$
for any nondecreasing convex function $f$, and hence for any function $f\in\F3$. \end{lemma}

\begin{lemma}\label{lem:incr in c}
If $X$ is a zero-mean r.v., then 
$\E f(c\,X)$ is nondecreasing in $c\ge0$ for any convex function $f$ and hence for any $f\in\F3$.   
\end{lemma}

\begin{lemma}\label{lem:p-monotone}
Let $\BS(p)\sim\bs(p)$. Then $e(p):=\E f(\BS(p))$ is nonincreasing in $p\in(0,1)$ for any $f\in\H2$ and so, 
by \eqref{eq:F-al-beta}, for any $f\in\H3$, whence, by Lemma~\ref{lem:F3=clG3}, for any $f\in\F3$.   
\end{lemma}

The extension from Theorem~\ref{th:bern} to Theorem~\ref{th:1} to 
Theorem~\ref{th:superm} is based in part on the following simple lemma, which may be of independent interest. 

\begin{lemma}\label{prop:superm} 
Suppose that for every $i\in\{1,\dots,n\}$ one has \eqref{eq:AiBi} and 
\eqref{eq:Bi/Ai}. 
Then 
\begin{align}
\E f(S_n) & \le\E f(c_1\BS_1+\dots+c_n\BS_n)
\label{eq:prop:superm}
\end{align}
for any $f\in\H2$, and so, 
by \eqref{eq:F-al-beta}, for any $f\in\H3$, whence, by Lemma~\ref{lem:F3=clG3} and Lebesgue's dominated convergence theorem, for any $f\in\F3$. 
\end{lemma}

For $m\ge1$, introduce
\begin{equation}\label{eq:p_*}
\begin{aligned}
p_*:=p_*(m)&:=\frac{2m+1 - \sqrt{4(m-1)(m+2)+1}}{4(2m-1)} \\
&= \frac{2}{(2m-1) \,(2m+1 + {\sqrt{4(m-1)(m+2)+1}}) },
\end{aligned}
\end{equation}
so that $p_*\in(0,\frac12]$.
Introduce also
\begin{equation*} 
\begin{aligned}
\de_1(u,c,p,m):=&\, 2c(1 - c^{2m-2})u+ 2pc(1 - c^{2m-1}) + c^2(1 - c^{2m-3}); \\
\de_2(u,c,p,m):=&\, 
(1-p)(1-{c^{2m-1}}){u^2} \\
& +2 c(1-{c^{2m-2}})u+2p c(1-{c^{2m-1}})+
{c^2}(1-{c^{2m-3}}); \\
\de_3(u,c,p,m):=&\, 
-{c^{2m-1}} {u^2}
   -2({c^{2m-1}}
    - c p+ {c^{2m}}
      p)\, u \\
& +\big((2 c +{c^2}
       -2 {c^{2m}}
        -{c^{2m+1}}) p-{c^{2m-1}}\big); \\
\de_4(u,c,p,m):=&\, (1-{c^{2m-1}})\, p\, {{(1+c+u)}^2}. 
\end{aligned}
\end{equation*}

\begin{lemma}\label{lem:equiv}
For any given pair $(p,m)$ such that $m\ge1$ and $p\in(0,1)$, the following two statements are equivalent to each other.
\emph{
\begin{enumerate}[(i)]
\item \emph{{\em (Cf. Statement \eqref{it:Bern-Sch,n=2} in Theorem~\ref{th:main}.) } 
For every every function 
$f\in\H3$, the function 
$$[0,\infty)^n\ni(a_1,a_2)\longmapsto
\E f(a_1^{1/(2m)}\BS_1+a_2^{1/(2m)}\BS_2),$$
is Schur-concave, where $\BS_i\iid\bs(p)$, $i=1,2$. }
\item
\emph{for all 
$u\in\R$ and $c\in(0,1)$, one has the inequalities }
\begin{align*}
& \de_1(u,c,p,m)\,\ii{u\ge0}\ge0; \\
& \de_2(u,c,p,m)\,\ii{-c\le u\le0}\ge0; \\
& \de_3(u,c,p,m)\,\ii{-1\le u\le-c}\ge0. 
\end{align*}
\end{enumerate}
}
\end{lemma}

\begin{lemma}\label{lem:u>0}
For all $u\ge0$, $c\in(0,1)$, $p\in[p_*,1)$, and $m\ge1$, one has 
$
\de_1(u):=\de_1(u,c,p,m)\ge0.
$ 
\end{lemma}

\begin{lemma}\label{lem:-c<u<0}
For all $u\in[-c,0]$, $c\in(0,1)$, $p\in[p_*,1)$, and $m\ge1$, one has 
$
\de_2(u):=\de_2(u,c,p,m)\ge0.
$
\end{lemma}

\begin{lemma}\label{lem:-1<u<-c}
For all $u\in[-1,-c]$, $c\in(0,1)$, $p\in[p_*,1)$, and $m\ge1$, one has 
$
\de_3(u):=\de_3(u,c,p,m)
\ge0.
$
\end{lemma}

\begin{lemma}\label{lem:p*,m*}
\emph{(Recall \eqref{eq:m_*} and \eqref{eq:p_*}.)}
For $p\in(0,1)$ and $m\ge1$, one has
$$m\ge m_*(p) \iff p\ge p_*(m). $$
\end{lemma}

\begin{lemma}\label{lem:IIIimpliesI}
In the context of Theorem~\ref{th:main}, implication 
$
\text{\eqref{it:Bern,n=2}}\implies\text{\eqref{it:m>m_*(p)}}
$
is true.
\end{lemma}

\begin{proof}[Proof of Theorem~\ref{th:main}]
It suffices to prove the implications
\begin{equation} \label{eq:implic}
\begin{gathered} 
\text{\eqref{it:m>m_*(p)}}\implies\text{\eqref{it:Bern-Sch,n=2}}\implies
\text{\eqref{it:Bern-Sch}}\implies\text{\eqref{it:Bern}}\implies
\text{\eqref{it:Bern,n=2}}\implies\text{\eqref{it:m>m_*(p)}}\quad\text{and} \\
\quad \text{\eqref{it:Bern}}\implies
\text{\eqref{it:superm}}\implies
\text{\eqref{it:bounded}}\implies\text{\eqref{it:Bern}}.
\end{gathered}
\end{equation} 

$\mathbf{
\text{\eqref{it:m>m_*(p)}}\implies\text{\eqref{it:Bern-Sch,n=2}}
}$:\quad 
Suppose that condition $m\ge m_*(p)$ of item \eqref{it:m>m_*(p)} takes place.
By Lemma~\ref{lem:p*,m*}, this condition is equivalent to $p\ge p_*(m)$. 
Now statement \eqref{it:Bern-Sch,n=2} with $\H3$ in place of $\F3$ follows from Lemmas \ref{lem:equiv}, \ref{lem:u>0}, \ref{lem:-c<u<0}, 
and \ref{lem:-1<u<-c}. 
Hence, by the definition \eqref{eq:GPlus3} of $\GPlus3$,
one has \eqref{it:Bern-Sch,n=2} with $\GPlus3$ in place of $\F3$.
To complete the proof of 
implication
$
\text{\eqref{it:m>m_*(p)}}\implies\text{\eqref{it:Bern-Sch,n=2}}
$,
it remains to use Lemma~\ref{lem:F3=clG3}; 
recall that the r.v.'s $\BS_i$ each take on only finitely many (namely, two) values.

$\mathbf{
\text{\eqref{it:Bern-Sch,n=2}}\implies\text{\eqref{it:Bern-Sch}}
}$:\quad 
By the well-known result by Muirhead \cite{muir} (see, e.g., \cite[Remark
B.1 of Chapter 2]{marshall}), a function of $n$ nonnegative arguments is Schur-concave iff it is Schur-concave in any two of its arguments. 
Now implication 
$\text{\eqref{it:Bern-Sch,n=2}}\implies\text{\eqref{it:Bern-Sch}}$ follows in view of Remark \ref{rem:shift} on page \pageref{rem:shift}, by conditioning on all of the r.v.'s $\BS_1,\dots,\BS_n$ except any given two of them.

$\mathbf{
\text{\eqref{it:Bern-Sch}}\implies\text{\eqref{it:Bern}}
}$:\quad 
Let here $a_1:=c_1^{2m},\dots,a_n:=c_n^{2m}$ and 
$\aa:=(a_1+\dots+a_n)/n$, so that 
$s^{(m)}=\aa\,{}^{1/(2m)}$. 
Note that $(a_1,\dots,a_n)\succeq(\underbrace{\aa,\dots,\aa}_n)$.
Now implication 
$
\text{\eqref{it:Bern-Sch}}\implies\text{\eqref{it:Bern}}
$ follows.

$\mathbf{
\text{\eqref{it:Bern}}\implies\text{\eqref{it:Bern,n=2}}
}$:\quad 
This implication is trivial.

$\mathbf{
\text{\eqref{it:Bern,n=2}}\implies\text{\eqref{it:m>m_*(p)}}
}$:\quad 
This implication is true by Lemma~\ref{lem:IIIimpliesI}.

$\mathbf{
\text{\eqref{it:Bern}}\implies\text{\eqref{it:superm}}
}$:\quad 
This implication follows immediately from Lemma~\ref{prop:superm}.

$\mathbf{
\text{\eqref{it:superm}}\implies\text{\eqref{it:bounded}}
}$:\quad 
This implication is trivial.

$\mathbf{
\text{\eqref{it:bounded}}\implies\text{\eqref{it:Bern}}
}$:\quad 
This implication is also trivial.
\end{proof}

\begin{proof}[Proof of Theorem~\ref{th:main-lefttail}]
Replace all functions $f$ in Theorem~\ref{th:main} with their reflections $\tilde f$ defined by $\tilde f(x):=f(-x)$ for all real $x$, replace all r.v.'s $X$ there with $-X$, and interchange $p$ with $q$, $a_i$ with $b_i$, and $A_{i-1}$ with $B_{i-1}$. Note that (i) $\tilde f(-X)=f(X)$; (ii) $\E X\le0\iff\E(-X)\ge0$; (iii) $X\sim\bs(p)\iff-X\sim\bs(q)$; and (iv) $(S_0,\dots,S_n)$ is a supermartingale with $S_0\le0$ a.s. iff $(-S_0,\dots,-S_n)$ is a submartingale with $-S_0\ge0$ a.s.
Thus,
Theorem~\ref{th:main-lefttail} follows immediately from Theorem~\ref{th:main}. 
\end{proof}

\begin{proof}[Proof of Proposition~\ref{prop:f3,ff3}]
This follows immediately from Lemmas~\ref{lem:f3,ff3} and \ref{lem:FF3neGG3}.
\end{proof}

\begin{proof}[Proof of Theorem~\ref{th:main-2tail}]
It suffices to prove the same implications, 
\eqref{eq:implic}, as in the proof of Theorem~\ref{th:main}, only with the changes stated in the formulation of Theorem~\ref{th:main-2tail}. 
Below, all these implications are understood in the context of Theorem~\ref{th:main-2tail}.
The proofs of most of these implications are similar to their proofs in the context of Theorem~\ref{th:main}. 
Below, only the most significant changes are described. 

$\mathbf{
\text{\eqref{it:m>m_*(p)}}\implies\text{\eqref{it:Bern-Sch,n=2}}
}$:\quad 
To prove this implication, in view of Theorems~\ref{th:main} and \ref{th:main-lefttail} and Proposition~\ref{prop:f3,ff3}, 
it suffices to verify that the function \eqref{eq:schur} is Schur-concave when $f(x)=x^2$ (for all real $x$) and $n=2$. Thus, it suffices to verify that, for any given 
$m\ge1$, the expression
$$\E(\BS_1\cos^{1/m}\th+\BS_2\sin^{1/m}\th)^2
=\cos^{2/m}\th+\sin^{2/m}\th$$
is nondecreasing in $\th\in[0,\pi/4]$.
But this is easy to see. 

$\mathbf{
\text{\eqref{it:Bern,n=2}}\implies\text{\eqref{it:m>m_*(p)}}
}$:\quad 
This implication follows from Theorems~\ref{th:main} and \ref{th:main-lefttail} and the observation that both classes $\F3$ and $\Fminus3$ are contained in $\FF3$. 

$\mathbf{
\text{\eqref{it:Bern}}\implies\text{\eqref{it:superm}}
}$:\quad  
In view of Theorems~\ref{th:main} and \ref{th:main-lefttail}, Proposition~\ref{prop:f3,ff3}, and Lebesgue's dominated convergence theorem, it suffices to verify that inequality \eqref{eq:Ef superm} holds when $f(x)=x^2$ for all real $x$ and $(S_0,\dots,S_n)$ is a martingale as described in the formulation of Theorem~\ref{th:main-2tail}. 
Note that inequality \eqref{eq:prop:superm} holds for the function $f_0(x):=x_+^2$ in place of $f$, since $f_0\in\G2$. It also holds for the function $\tilde f_0(x):=(-x)_+^2$ in place of $f$, given that $(S_0,\dots,S_n)$ is a martingale as described. Thus, \eqref{eq:prop:superm} holds when $f(x)=x^2$ for all real $x$. It remains to note that    
$$
\E(c_1\BS_1+\dots+c_n\BS_n)^2=n\,(s^{(1)})^2\le n\,(s^{(m)})^2
=\E(s^{(m)}\cdot(\BS_1+\dots+\BS_n))^2,
$$
so that one does have inequality \eqref{eq:Ef superm} when $f(x)=x^2$ for all real $x$.
\end{proof}

\subsection{Proofs of the lemmas}\label{proofs2} 

\begin{proof}[Proof of Lemma~\ref{lem:g',g'',g'''}]
The convexity of $g$ implies $g(0)-g(x)\ge g'(x)(-x)$ and 
$g(2x)-g(x)\ge g'(x)x$, so that $|g'(x)||x|\le\max(|g(0)-g(x)|,|g(2x)-g(x)|)$ for all $x\in\R$.
Letting now 
$x\to-\infty$ and using the existence of the finite limit $g(-\infty)$, one has $g'(-\infty)=0$. 
\end{proof}

\begin{proof}[Proof of Lemma~\ref{lem:F-al}]
This lemma was stated essentially as Proposition~1.1 in \cite{binom}. The proof given here is a little more detailed.  
Assume first that $f\in\H\al$, so that
$f(x)=  
\int (x-t)_+^\al\,\d\mu(t)$ for a function $\mu\in\M\al$, whence 
$f^{(\al-1)}(x)=
\al!\int (x-t)_+\,\d\mu(t)$ is convex as a limit of linear combinations with nonnegative coefficients of convex functions 
$x\mapsto(x-t)_+$. The conditions 
$f^{(0)}(-\infty)=\dots=f^{(\al-1)}(-\infty)=0$ follow by Lebesgue's dominated convergence theorem. Thus, $f\in\tilde{\H\al}$.
Moreover, it is clear that all the functions $f^{(0)},\dots,f^{(\al-1)}$ are nonnegative. 
 
Assume now that $f\in\tilde{\H\al}$. 
Consider first the case $\al=1$. Then $f$ is convex and $f(-\infty)=0$. Hence, by Lemma~\ref{lem:g',g'',g'''}, one has $f'(-\infty)=0$. Therefore, 
$$f(x)=\smallint_{-\infty}^x f'(u)\,\d u=\smallint_{-\infty}^x \d u\,\smallint_{-\infty}^u \d\,f'(v)
=\smallint (x-v)_+\,\d\,\mu(v),$$
by Fubini's theorem, where $\mu:=f'$; thus, $f\in\H1$. 
The case of any natural $\al\ge2$ can now be treated by induction, in a similar manner. 
Indeed, if $f\in\tilde{\H\al}$ for a natural $\al\ge2$, then 
$f'\in\tilde{\mathcal{H}}_+^{\al-1}$, by the definition of $\tilde{\H\al}$. Hence, for a function $\mu\in\M\al$,
$$f(x)=\smallint_{-\infty}^x f'(u)\,\d u
=\smallint_{-\infty}^x \d u\,\smallint (u-v)_+^{\al-1}\,\d\mu(v)
=\smallint (x-v)_+^\al\,\d\,\mu(v)/\al,$$
so that $f\in\H\al$.
\end{proof}

\begin{proof}[Proof of Lemma~\ref{lem:f-in-clG3}]
For the given function $f$ and any $y\in\R$, introduce the functions defined by the formulas
\begin{align}
& f_{2,y}(x) :=\big(f''(y)+f'''(y)(x-y)\big)_+\, \ii{x\le y} + f''(x)\,\ii{x>y}; 
\label{eq:f_{2,y}}\\
&
\begin{aligned}
f_y(x) :=
\Big(f(y)+f'(y)(x-y)+\smallint_{-\infty}^y f_{2,y}(u)\,(u-x)_+\,\d u
\Big)\,& \ii{x\le y} \\
+ f(x)\,& \ii{x>y}
\end{aligned}
\label{eq:f_y}
\end{align}
for all real $x$.
Here $f'''$ denotes the right derivative of the convex function $f''$, so that $f'''$ is nondecreasing. 
Note that
\begin{equation} \label{eq:f''_y=f_{2,y}}
f''_y=f_{2,y}.
\end{equation}

Since $f''$ is convex, one has 
\begin{equation} \label{eq:f''over tangent}
f''(x)\ge f''(y)+f'''(y)(x-y)
\end{equation} 
for all $y$ and $x$; also, it is given that $f''$ is nonnegative; it follows that 
\begin{equation} \label{eq:f''>f_{2,y}}
f''\ge f_{2,y}.
\end{equation} 

Observe that, moreover, the family of functions $(f_{2,y})$ is nonincreasing in 
$y\in\R$. 
Indeed, let $y$ and $y_1$ be any real numbers such that $y_1<y$. 
Then $f_{2,y_1}=f''\ge f_{2,y}$ on $[y_1,\infty)$, in view of \eqref{eq:f_{2,y}} and 
\eqref{eq:f''>f_{2,y}}. 
Recalling \eqref{eq:f''over tangent} and the fact that $f'''$ is nondecreasing, one has the inequalities $f''(y_1)\ge f''(y)+f'''(y)(y_1-y)$ and 
$f'''(y_1)(x-y_1)\ge f'''(y)(x-y_1)$ for all $x\le y_1$; adding these inequalities, one sees that 
$f''(y_1)+f'''(y_1)(x-y_1)\ge f''(y)+f'''(y)(x-y)$. 
It follows, in view of \eqref{eq:f_{2,y}}, that $f_{2,y_1}\ge f_{2,y}$ on the interval $(-\infty,y_1]$ as well, and hence on the entire real line.   

Using integration-by-parts/Fubini's theorem as in the proof of Lemma~\ref{lem:F-al}, one can verify that for any function $g\in\C^2$ and all real $y$ and $x$
\begin{equation} \label{eq:g,g''}
\begin{aligned}g(x)=\Big(g(y)+g'(y)(x-y)+\smallint_{-\infty}^y g''(u)\,(u-x)_+\,\d u
\Big)\,& \ii{x\le y} \\
+ g(x)\,& \ii{x>y}.
\end{aligned}
\end{equation}
By \eqref{eq:f''_y=f_{2,y}}, for any real $w$ one has $f_w\in\C^2$, so that one can substitute $f_w$ for $g$ in \eqref{eq:g,g''}.
In fact, let us do so for $w\in\{y,y_1\}$, again assuming that $y_1<y$. 
At that, by \eqref{eq:f_y}, one has $f_{y_1}=f=f_y$ on the interval $[y,\infty)$ and hence $f_{y_1}(y)=f_y(y)$ and $f'_{y_1}(y)=f'_y(y)$.
Now, since the family of functions $(f''_y)_{y\in\R}=(f_{2,y})_{y\in\R}$ is nonincreasing, one can see that the family $(f_y)_{y\in\R}$ is nonincreasing as well. 
Next, since $\ii{x\le y}\to0$ and $\ii{x>y}\to1$ for each $x$ as $y\to-\infty$, one concludes, in view of  \eqref{eq:f''_y=f_{2,y}},
that for any decreasing sequence $(y_n)$ in $\R$ converging to $-\infty$ one has 
$f_{y_n}\to f$, in the sense of Definition~\ref{def:convergence}. 

It remains to verify that for every real $y$ one has $f_y\in\G3$ and, moreover, 
$f_y\in\GPlus3$ in the case when $f$ is known to be nondecreasing. Observe that
\begin{equation} \label{eq:f_{2,y}cases}
f_{2,y}(x)=f'''(y)\,(x-z)\,\ii{z\le x\le y}+f''(x)\,\ii{x>y},
\end{equation}
where 
$$z:=y\,\ii{f'''(y)=0}+(y-f''(y)/f'''(y))\,\ii{f'''(y)\ne0}.$$
Indeed, $f'''$ is nonnegative and nondecreasing (since $f''$ is nondecreasing and convex). 
Hence, in the case when $f'''(y)=0$, one has $f'''=0$ on the entire interval 
$(-\infty,y]$. This and the condition $f''(-\infty)=0$ implies $f''=0$ on the entire interval $(-\infty,y]$, so that $f''(y)=0$. Now one sees that expressions 
\eqref{eq:f_{2,y}} and \eqref{eq:f_{2,y}cases} both equal $f''(x)\,\ii{x>y}$ in the case when $f'''(y)=0$.
In the other case, when $f'''(y)\ne0$, one has $f'''(y)>0$ 
(since $f'''$ is nonnegative). Also, here 
$f''(y)+f'''(y)(x-y)=f'''(y)\,(x-z)$, whence \eqref{eq:f_{2,y}cases} again follows.

Now, for the right derivative 
$f'_{2,y}$ of $f_{2,y}$, \eqref{eq:f_{2,y}cases} yields
$$f'_{2,y}(x)=f'''(y)\,\ii{z\le x\le y}+f'''(x)\,\ii{x>y}.$$
Since $f'''$ is nonnegative and nondecreasing, it follows now that $f'_{2,y}$ is nondecreasing. 
Therefore, $f_{2,y}$ is convex. 
That is, by \eqref{eq:f''_y=f_{2,y}}, $f''_y$ is convex. 
Also, \eqref{eq:f_{2,y}cases} and \eqref{eq:f''_y=f_{2,y}}
show that $f''_y=0$ on the interval 
$(-\infty,z]$. This means that 
\begin{equation} \label{eq:f_y=a+bx}
\text{$f_y(x)=a+b\,x$ for some real constants $a$ and $b$ and all $x\le z$.}
\end{equation} 

Let now 
$$h_y(x):=f_y(x)-(a+b\,x)$$
for all real $x$. Then $h''_y=f''_y$ is convex. Moreover, $h_y=0$ on the interval 
$(-\infty,z]$, so that $h_y(-\infty)=h'_y(-\infty)=h''_y(-\infty)=0$. 
By Lemma~\ref{lem:F-al}, $h_y\in\H3$. 
Thus, $f_y\in\G3$. 

If, moreover, $f$ is nondecreasing, then $f'\ge0$. Hence and because (in view of 
\eqref{eq:f''>f_{2,y}})
$f''(u)\ge f_{2,y}(u)\ge f_{2,y}(u)\,\ii{u>x}$ for all $x$, $y$, and $u$, one has
\begin{align*}
f'_y(x) & =
\Big(f'(y)-\smallint_{-\infty}^y f_{2,y}(u)\,\ii{u>x}\,\d u\Big)\,\ii{x\le y}
			+f'(x)\,\ii{x>y} \\
& \ge f'(-\infty)\,\ii{x\le y}
			+f'(x)\,\ii{x>y} \ge0			
\end{align*}
for all $x$. Now \eqref{eq:f_y=a+bx} implies $b\ge0$. 
Since $h_y\in\H3$, one finally sees that 
$f_y\in\GPlus3$.
\end{proof}

\begin{proof}[Proof of Lemma~\ref{lem:F3=clG3}]
First note that $\cl\F3=\F3$, because the pointwise convergence preserves both the monotonicity and the convexity. 

Next, take any $f\in\GPlus3$, so that 
$$f(x)=a+b\,x+\smallint(x-t)_+^3\,\d\mu(t)$$
for all $x$, where $a\in\R$, $b\ge0$, and $\mu$ is nondecreasing and 
$\smallint(-t)_+^3\,\d\mu(t)<\infty$. 
It follows that $f$ is nondecreasing and convex, since the functions 
$x\mapsto(x-t)_+^3$ are so. Similarly, $f''$ is nondecreasing and convex, since 
$f''(x)=6\,\smallint(x-t)_+\,\d\mu(t)$.
That is, $f\in\F3$ for any $f\in\GPlus3$, so that $\GPlus3\subseteq\F3$, whence 
$\cl\GPlus3\subseteq\cl\F3=\F3$. 
 
It remains to show that $\F3\subseteq\cl\GPlus3$. Take any $f\in\F3$. Then, by definition \eqref{eq:F3}, $f$ and $f''$ are nondecreasing and convex. Hence, $f'$ is nonnegative, nondecreasing, and convex. Now Lemma~\ref{lem:g',g'',g'''} yields 
$f''(-\infty)=0$. 
Also, $f''$ is nonnegative, since $f$ is convex.
Thus, by Lemma~\ref{lem:f-in-clG3}, $f\in\cl\GPlus3$. 
\end{proof}

\begin{proof}[Proof of Lemma~\ref{lem:f3,ff3}]
First note that $\cl\FF3=\FF3$, because the pointwise convergence preserves the convexity. 

Next, it is trivial that $\GG3\subseteq\FF3$, whence $\cl\GG3\subseteq\cl\FF3=\FF3$.

It remains to show that $\FF3\subseteq\cl\GG3$. Take any $f\in\FF3$. Then, by definition \eqref{eq:FF3}, $f$ and $f''$ are convex. The latter condition implies that at least 
one of the following three cases must take place: $f''$ is nondecreasing on $\R$ or 
$f''$ is nonincreasing on $\R$ or $f''$ switches from nonincreasing to nondecreasing. 

{\em Case 1: $f''$ is nondecreasing on $\R$.} \quad
Since $f$ is convex, $f''\ge0$ on $\R$. Hence, there exists the limit
$c:=f''(-\infty)\in[0,\infty)$. Let 
$$g(x):=f(x)-c\,x^2/2$$
for all real $x$. Then $g''=f''-c=f''-f''(-\infty)\ge0$, since $f''$ is nondecreasing. 
Also, $g''=f''-c$ is nondecreasing and convex, since $f''$ is so. In addition, $g''(-\infty)=0$. 
Therefore, by Lemma \ref{lem:f-in-clG3}, $g\in\cl\G3$. That is, there exists a sequence of functions $(g_n)$ such that $g_n\to g$ and
$$g_n(x)=a_n+b_n\,x+h_n(x)$$
for all real $x$, where, for each $n$, $a_n$ and $b_n$ are real constants and $h_n$ is a function in $\H3$; then $h_n\in\F3$ (because $\H3\subseteq\F3$, as seen, for example, from the proof of Lemma~\ref{lem:F-al}). 
Let now
$$f_n(x):=c\,x^2/2+g_n(x)
=c\,x^2/2+a_n+b_n\,x+h_n(x)$$
for all $x$ and $n$.
Then $g_n\to g$ implies $f_n\to f$. Moreover, for every $n$ one has $f_n\in\GG3$. 
Indeed, if $b_n<0$, then the function $h_n$ belongs to $\F3$ and the function
$x\mapsto a_n+b_n\,x$ belongs to $\Fminus3$; and if $b_n\ge0$, then the function 
$x\mapsto a_n+b_n\,x+h_n(x)$ belongs to $\F3$ and the function 
$x\mapsto0$ belongs to $\Fminus3$.
Thus, $f\in\cl\GG3$ for any $f\in\FF3$ satisfying the condition of Case~1. 

{\em Case 2: $f''$ is nonincreasing on $\R$.} \quad This case reduces to Case~1 by considering the function $x\mapsto \tilde f(x):=f(-x)$ in place of $f$. Indeed, if $f\in\FF3$ and $f''$ is nonincreasing, then 
$\tilde f\in\FF3$ and $\tilde f''$ is nondecreasing. Moreover, 
$f_n\in\GG3\iff\tilde f_n\in\GG3$, where $\tilde f_n(x):=f_n(-x)$ for all real $x$. 

{\em Case 3: There exists some real $x_0$ such that $f''$ is nonincreasing on $(-\infty,x_0]$ and nondecreasing on $[x_0,\infty)$.} \quad Here without loss of generality (w.l.o.g.) $x_0=0$. 
Let
$$g(x):=f(x)-f'(0)x-f''(0)x^2/2$$
for all real $x$, so that $g(0)=f(0)$ and $g'(0)=g''(0)=0$; moreover, $g''=f''-f''(0)$ is 
convex on $\R$ (since $f''$ is so),
nonincreasing on $(-\infty,0]$, and nondecreasing on $[0,\infty)$, whence $g''\ge0$. 
Let, for all real $x$,
$$h_+(x):=g''(x)\ii{x>0}\quad\text{and}\quad
h_-(x):=g''(x)\ii{x\le0},$$
so that $h_+ + h_- =g''$, $h_\pm\ge0$, $h_+$ is nondecreasing, and $h_-$ is nonincreasing. 
Also, $h_+$ and $h_-$ are convex, since $g''$ is convex and $g''\ge g''(0)=0$. 
Let further, for all real $x$,
$$H_+(x):=\smallint\,(x-t)_+\,h_+(t)\,\d t\quad\text{and}\quad
H_-(x):=\smallint\,(t-x)_+\,h_-(t)\,\d t,$$
so that
$H_\pm(0)=0$, 
\begin{align*}
H'_+(x) & =\smallint\,\ii{x>t}\,h_+(t)\,\d t
=\smallint_0^x\,g''(t)\,\d t\,\ii{x>0}
=g'(x)\ii{x>0}; \\
H''_+(x) & =g''(x)\ii{x>0}=h_+(x);\\
H'_-(x) & =-\smallint\,\ii{t>x}\,h_-(t)\,\d t
=-\smallint_x^0\,g''(t)\,\d t\,\ii{x\le0}
=g'(x)\ii{x\le0}; \\
H''_-(x) & =g''(x)\ii{x\le0}=h_-(x).
\end{align*}
It follows that $H'_+ + H'_-=g'$, whence for all real $x$ one has
$H_+(x)+H_-(x)=g(x)-g(0)$, that is,
\begin{equation} \label{eq:f(x)=}
f(x)=f(0)+f'(0)x+f''(0)x^2/2+H_+(x)+H_-(x).
\end{equation}
Also, $H_+''=h_+$ is nonnegative, nondecreasing, and convex, and hence $H_+$ is also convex. Also, since $h_+\ge0$, the first expression for $H'_+(x)$ above shows that $H'_+\ge0$.
Thus, $H_+$ and $H_+''$ are nondecreasing and convex; that is, $H_+\in\F3$. 
Similarly, $H_-\in\Fminus3$.  

Note also that $f''(0)\ge0$, since $f$ belongs to $\FF3$ and is hence convex. 
If $f'(0)<0$, then the function $H_+$ belongs to $\F3$ and the function 
$x\mapsto f(0)+f'(0)x+H_-(x)$ belongs to $\Fminus3$; and if $f'(0)\ge0$, then the function $x\mapsto f(0)+f'(0)x+H_+(x)$ belongs to $\F3$ and the function
$H_-$ belongs to $\Fminus3$.
Thus, \eqref{eq:f(x)=} implies that
$f\in\GG3\subseteq\cl\GG3$ for any $f\in\FF3$ satisfying the condition of Case~3. 

One concludes that, in all cases $f\in\FF3$ implies $f\in\cl\GG3$. That is, 
$\FF3\subseteq\cl\GG3$ indeed. 
\end{proof}

\begin{proof}[Proof of Lemma~\ref{lem:f at infty}]
Let $f\in\F3$. 

{\em Case 1: $f''=0$ on $\R$.}\quad Then there exist real $a$ and $b$ such that $f(x)=a+b\,x$ for all real $x$, so that $f(x)=O(x)$ as $x\to\infty$. 

{\em Case 2: there exists some $t\in\R$ such that $f''(t)\ne0$.}\quad By 
\eqref{eq:F3}, 
$f''$ is nonnegative (because $f$ is convex) and nondecreasing. Hence, $f''(t)>0$ and $f''(x)\ge f''(t)$ for all $x\ge t$. It follows that
$f(x)\ge f(t)+f'(t)\,(x-t)+f''(t)\,(x-t)^2/2$ for all $x\ge t$, whence 
$\liminf_{x\to\infty}f(x)/x^2\ge f''(t)/2>0$.
\end{proof}

\begin{proof}[Proof of Lemma~\ref{lem:f at -infty}]
Let $f\in\F3$. Since $f$ is nondecreasing, one has $f(x)\le f(0)$ for all $x\le0$. On the other hand, $f(x)\ge f(0)+f'(0)x$ for all real $x$, since $f$ convex. It follows that $|f(x)|\le|f(0)|+|f'(0)x|$ for all $x\le0$, so that $f(x)=O(|x|)$ as 
$x\to-\infty$.
\end{proof}

\begin{proof}[Proof of Lemma~\ref{lem:FF3neGG3}]
For all real $x$, let
\begin{equation} \label{eq:f example}
f(x):=\tfrac83\,(1-x)^{3/2}\,\ii{x\le0}
+(\tfrac83-4x+x^2+\tfrac16\,x^3+\tfrac1{16}\,x^4)\,\ii{x>0}.
\end{equation}
Then it is easy to see that $f\in\FF3$. 

Suppose that $f\in\GG3$. Then, by \eqref{eq:GG3}, there exist $c\ge0$, $f_+\in\F3$, and $f_-\in\Fminus3$ such that for all real $x$
$$f(x)=c\,x^2/2+f_+(x)+f_-(x).$$
Let $x\to-\infty$. Then,
by Lemma~\ref{lem:f at -infty}, $f_+(x)=O(|x|)$, while by 
Lemma~\ref{lem:f at infty}, 
either $f_-(x)=O(|x|)$ or 
$\liminf_{x\to-\infty}f_-(x)/x^2>0$. 
It follows that either $f(x)=O(|x|)$ as $x\to-\infty$ or 
$\liminf_{x\to-\infty}f(x)/x^2>0$. However, neither of these two alternatives is compatible with the fact that $f(x)=\tfrac83\,(1-x)^{3/2}$ for $x\le0$. 
Thus, $f\in\FF3\setminus\GG3$.
\end{proof}

\begin{proof}[Proof of Proposition~\ref{prop:F3neG3}]
Let $g:=f'$, where $f$ is defined by \eqref{eq:f example}. Then it is easy to see that $g\in\F3$. 

Suppose that $g\in\G3$. Then, by \eqref{eq:G3}, there exist real $a$ and $b$ and $h\in\H3$ such that for all real $x$ one has $g(x)=a+b\,x+h(x)$ and hence $g'(x)=b+h'(x)$. By 
Lemma~\ref{lem:F-al}, $h(-\infty)=h'(-\infty)=0$. Hence, $b=g'(-\infty)$. But, by inspection, $g'(-\infty)=0$. It follows that $b=0$, and so, 
$g(-\infty)=a+h(-\infty)=a\in\R$. This contradicts the fact that 
$g(-\infty)=-\infty$.  
\end{proof}

\begin{proof}[Proof of Lemma~\ref{lem:extreme-bounded}]
Write
$$X=\frac{X+a}{b+a}\,b+\frac{b-X}{b+a}\,(-a).$$
Let $f$ be any nondecreasing convex function. The convexity (together with the condition $-a\le X\le b$ a.s.) implies that 
$$f(X)\le\frac{X+a}{b+a}\,f(b)+\frac{b-X}{b+a}\,f(-a)
=\E f(\sqrt{a\,b}\,\BS)+\frac{f(b)-f(-a)}{b+a}\,X$$
a.s. Now, since $f$ is nondecreasing and $\E X\le0$, the lemma follows. 
\end{proof}

\begin{proof}[Proof of Lemma~\ref{lem:incr in c}]
Since $f$ is convex, the function $[0,\infty)\ni c\mapsto f(c\,X)$ is convex as well. Hence, the function
$[0,\infty)\ni c\mapsto g(c):=\E f(c\,X)$ is convex. Since $\E X=0$, one has 
$g(c)\ge g(0)$ for all real $c$, by Jensen's inequality. Therefore, the right derivative of $g$ is nonnegative at $0$ and hence on $[0,\infty)$. Now the lemma follows. 
\end{proof}

\begin{proof}[Proof of Lemma~\ref{lem:p-monotone}]
In view of the definition of the class $\H2$, it suffices to verify the statement of the lemma for all functions $f$ of the form $f_t(x):=(x-t)_+^2$, for all real $t$, so that $e(p)=\E(\BS(p)-t)_+^2$. 
Then 
\begin{align*}
e'(p) & =\lp\sqrt{\tfrac qp}-t\rp\, 
\lp-\sqrt{\tfrac pq}-t\rp\,\ii{-\sqrt{\tfrac pq}<t<\sqrt{\tfrac qp}} \le0.
\end{align*} 
\end{proof}

\begin{proof}[Proof of Lemma~\ref{prop:superm}]
W.l.o.g., the $\BS_i$'s are independent of the $S_i$'s. 
For $i=0,1,\dots,n$ and $f\in\H2$, introduce
$$F_i:=\E f\lp S_i+c_{i+1}\BS_{i+1}+\dots+c_n\BS_n\rp.$$
Recall that, by Remark \ref{rem:shift}, the classes $\H\al$ are invariant with respect to the shifts. Hence, by 
Lemmas~\ref{lem:extreme-bounded}, 
\ref{lem:incr in c}, and \ref{lem:p-monotone},
\begin{align*}
F_i & = \E\E_{i-1} f(S_{i-1}+X_i
+c_{i+1}\BS_{i+1}+\dots+c_n\BS_n) \\
& \le \E\E_{i-1} f(S_{i-1}+\sqrt{A_{i-1}B_{i-1}}\,\tilde{\BS_i}
+c_{i+1}\BS_{i+1}+\dots+c_n\BS_n) \\
& \le \E\E_{i-1} f(S_{i-1}+c_i\,\tilde{\BS_i}
+c_{i+1}\BS_{i+1}+\dots+c_n\BS_n) \\
& \le \E\E_{i-1} f(S_{i-1}+c_i\BS_i
+c_{i+1}\BS_{i+1}+\dots+c_n\BS_n)\\
& = F_{i-1}
\end{align*} 
for $i=1,\dots,n$,
where $\E_{i-1}$ denotes the conditional expectation given the $\si$-algebra 
$G_{i-1}$ generated by $H_{i-1}$ and 
$(\BS_{i+1},\dots,\BS_n)$, and 
the conditional distribution of $\tilde{\BS_i}$ given $G_{i-1}$ is $\bs(p_i)$, with 
$$p_i:=\frac{A_{i-1}}{B_{i-1}+A_{i-1}},$$
so that $p_i\ge p$, according to \eqref{eq:Bi/Ai}. Hence,
$$
\E f(S_n)=F_n\le F_0\le\E f(c_1\BS_1+\dots+c_n\BS_n);
$$
the last inequality follows because $S_0\le0$ a.s. and all functions $f\in\H2$ are nondecreasing. 
\end{proof}

\begin{proof}[Proof of Lemma~\ref{lem:equiv}]
Statement (i) is equivalent to the following: 
for every $t\in\R$, the function 
$$[0,\infty)^2\ni(a_1,a_2)\longmapsto
e_{t,p,m}(a_1,a_2):=\E f_t(a_1^{1/(2m)}\BS_1+a_2^{1/(2m)}\BS_2)$$
is Schur-concave, where $\BS_i\iid\bs(p)$ and
$$f_t(x):=\tfrac13\,(x-t)_+^3.$$
Using the homogeneity property 
$$e_{t,p,m}(\la a_1,\la a_2)=\la^{3/(2m)}\,e_{\la^{-1/(2m)}\,t,p,m}(a_1,a_2)$$
for every $\la>0$, one may assume w.l.o.g. that 
$a_1+a_2=1$, so that $a_1=\cos^2\th$ and $a_2=\sin^2\th$, for some 
$\th\in[0,\pi/2]$; moreover, in view of the same homogeneity property,  one may replace here the i.i.d. standardized Bernoulli r.v.'s $\BS_1$ and $\BS_2$ with i.i.d. {\em centered} Bernoulli r.v.'s $\BC_1$ and $\BC_2$, such that 
$$\P(\BC_i=1-p)=p=1-\P(\BC_i=-p),\quad i=1,2.$$
Therefore, statement (i) is equivalent to 
\begin{equation} \label{eq:de}
\De_{p,m}(\th,t)
:=\partial_\th\,\big(\E f_t(\BC_1\,\cos^{1/m}\th+\BC_2\,\sin^{1/m}\th)\big)
\end{equation}
being nonnegative for all $\th\in[0,\pi/4]$ and all $t\in\R$ (where 
$\partial_\th:=\partial/\partial\th$), which is in turn 
equivalent to 
\begin{multline} \label{eq:De1,De2}
\De_{1,p,m}(\th,u):=
\frac{m}{(1-p)p}\,\cos^{1-1/m}\th\,\sin^{1-1/m}\th\, \\
\times\De_{p,m}\big(\th,-u-p\,\cos^{1/m}\th-p\,\sin^{1/m}\th\big) \\
=\De_{2,p,m}(\cos^{1/m}\th,\sin^{1/m}\th,u)
\end{multline}
being nonnegative for all $\th\in[0,\pi/4]$ and $u\in\R$, where 
\begin{multline*}
\De_{2,p,m}(c_1,c_2,u):=\\
-({c_1^{2m-1}}
  -{c_2^{2m-1}}
   )\,q {\,u_+^2}
    -({c_2^{2m-1}}
      q+{c_1^{2m-1}}
       p) {\,(c_1+u)_+^2} \\
       +({c_1^{2m-1}}
         q+{c_2^{2m-1}}
          p) {\,(c_2+u)_+^2}
          +({c_1^{2m-1}}
           -{c_2^{2m-1}}
            )\, p {\,(c_1+c_2+u)_+^2
               }.
\end{multline*}
Now, in view of the homogeneity relation
$$\De_{2,p,m}(c_1,c_2,u)=c_1^{2m+1}\,\De_{2,p,m}(1,c_2/c_1,u/c_1),
\quad\text{where}\quad 0<c_2<c_1,$$ 
statement (i) reduces to 
$\De_{2,p,m}(1,c,u)$ being nonnegative for all $c\in(0,1)$ and $u\in\R$.

It remains to note that
\begin{equation}\label{eq:cases}
\De_{2,p,m}(1,c,u)=
\begin{cases}
\de_1(u,c,p,m)\quad&\text{if}\quad u\ge0, \\
\de_2(u,c,p,m)\quad&\text{if}\quad u\in[-c,0], \\
\de_3(u,c,p,m)\quad&\text{if}\quad u\in[-1,-c], \\
\de_4(u,c,p,m)\quad&\text{if}\quad u\in[-1-c,-1], \\
0\quad&\text{otherwise},
\end{cases}
\end{equation}
and
$\de_4(u,c,p,m)$
is manifestly nonnegative for all 
$c\in(0,1)$, $u\in\R$, $p\in(0,1)$, and $m>1$.
\end{proof}

\begin{proof}[Proof of Lemma~\ref{lem:u>0}]
Note that $\de_1(u)\ge\de_1(0)$ for $u\ge0$, $c\in(0,1)$, and $m\ge1$.
Next, \eqref{eq:cases} shows that $\de_1(0)=\de_2(0)$.
Now the lemma follows immediately from Lemma~\ref{lem:-c<u<0}.
\end{proof}

\begin{proof}[Proof of Lemma~\ref{lem:-c<u<0}]
W.l.o.g., $m>1$. Note that 
$$\partial_p\de_2(u,c,p,m)=(2c-u^2)(1-{c^{2m-1}})>0$$ 
for $u\in[-c,0]$, $c\in(0,1)$, and $m>1$,
so that w.l.o.g. 
$$p=p_*=p_*(m).$$
Next,
$\de_2(u)$ is a convex quadratic polynomial, whose minimum over all $u\in\R$ is attained at
$$u=u_*(c,p,m):=-\frac{c(1-{c^{2m-2}})}
    {(1-{c^{2m-1}})(1-p)}.$$
Hence, it suffices to show that $\g(c,p_*)$ is nonnegative for all $c\in[0,1]$ and $m>1$, where
\begin{multline*}
\g(c,p):=\de_2(u_*(c,p,m),c,p,m)\,\frac{(1-{c^{2m-1}})(1-p)}{c} \\
=-2 {{\big(1-{c^{2m-1}}\big)}^2}
   {p^2}+(1-{c^{2m-1}}
    ) (2 -c+{c^{2m-2}}
     -2 {c^{2m-1}}
      ) p-{c^{2m-2}}{{(1-c)}^2}
\end{multline*}
The main idea in the proof of this lemma is to replace here the entry of $p^2$ with the equivalent (for $p=p_*$), first-degree in $p$ polynomial expression according to the identity
\begin{equation}\label{eq:p^2}
p_*^2=
\frac{(4{m^2}-1)p_*-1}
   {2{{(2m-1)}^2}},
\end{equation}
to obtain
\begin{equation*} 
(2m-1)^2 \g(c,p_*)=f(c,p_*),
\end{equation*}
where 
\begin{multline*}
f(c,p):=
(1-{c^{2m-1}}{)^2}
   -{{(2m-1)}^2}
    {{(1-c)}^2} {c^{2m-2}} \\
       +p\,(1-{c^{2m-1}}
       ) (2 m-1) \big((2m-3)(1- {c^{2m-1}}
         )-(2m-1)\,c\,(1- {c^{2m-3}})\big).
         \end{multline*}
It suffices to show that $f(c,p)\ge0$ for all $p\in(0,1)$, $m>1$, and 
$c\in(0,1)$.
         Introduce
\begin{align*}
g(c)&:=f(c,p)/c^{2m}; \\
g_1(c)&:=g'(c)/c^{2m-3}; \\
g_2(c)&:=g_1'(c)\,c^2; \\
g_3(c)&:=g_2'(c)/c^{1-2m}; \\
g_4(c)&:=g_3'(c)/c^{1-2m}; \\
g_5(c)&:=g_4'(c)/c^{2m-4}.
\end{align*}
Then, letting
$$s:=m-1>0,$$ 
one has
$g_5'(c)=
8{c^{-1-2s}}\,s\,(1+s)(1+2s{)^2}(1+4s)\big(1-p+4{s^2}p\big)>0$ for all $c\in(0,1)$, and so, $g_5$ is increasing on $(0,1)$ to 
$$g_5(1)=-16\,s\,(1+s)(1+2s)^2\big(1-p+s+4{s^2}p\big)<0.$$  
Hence, $g_5<0$ on $(0,1)$, so that $g_4$ is decreasing on $(0,1)$ to 
$$g_4(1)=8\,s\,(1+2s)^2\big(1-p+s+4{s^2}p\big)>0.$$
Hence, $g_4>0$ on $(0,1)$. Since $g_3(1)=g_2(1)=g_1(1)=g(1)=0$, it follows successively that $g_3<0$, $g_2>0$, $g_1<0$, and $g>0$ on $(0,1)$.
This completes the proof of 
Lemma~\ref{lem:-c<u<0}. 
\end{proof}

\begin{proof}[Proof of Lemma~\ref{lem:-1<u<-c}]
This follows because $\de_3(u)$ is concave in $u$, 
$\de_3(-1)={c^2}(1-{c^{2m-1}}) p\ge0$, and 
$\de_3(-c)=\de_2(-c)\ge0$,
where the latter equality and inequality follow immediately from 
\eqref{eq:cases} and Lemma~\ref{lem:-c<u<0}, respectively.
\end{proof}

\begin{proof}[Proof of Lemma~\ref{lem:p*,m*}]
It is clear from the second expression for $p_*(m)$ in \eqref{eq:p_*} that 
$p_*(m)$ decreases continuously from $\frac12$ to $0$ as $m$ increases from $1$ to $\infty$. Also, one can verify that $m_*(p_*(m))=m$ for all $m\ge1$ (here one may use identity \eqref{eq:p^2}). 

If now $p\in(0,p_*(m))$, then $p=p_*(m_1)$ for some $m_1>m$, whence 
$m<m_1=m_*(p_*(m_1))=m_*(p)$.

It remains to consider the condition $p\ge p_*(m)$. If at that $p>\frac12$, then 
$m\ge1=m_*(p)$, by \eqref{eq:m_*}. 
If, however, $p\in[p_*(m),\frac12]$, then $p=p_*(m_1)$ for some 
$m_1\in[1,m]$, whence $m\ge m_1=m_*(p_*(m_1))=m_*(p)$.
\end{proof}

\begin{proof}[Proof of Lemma~\ref{lem:IIIimpliesI}]
Suppose, to the contrary, that statement \eqref{it:Bern,n=2} of 
Theorem~\ref{th:main} is true, while $m<m_*(p)$. Then, by 
Lemma~\ref{lem:p*,m*}, one has
$$p<p_*:=p_*(m);$$
then, in particular, one has $0<p<\frac12$.  
Introduce 
$$
u_p:=-\frac{{2^{1-\frac{1}{2m}}}
      (m-1)}{(2m-1)(1-p)}.
$$
In view of the elementary inequality $p_*\le1/(2m-1)$ for $m\ge1$ and the condition $p<p_*$, one has 
$$p<\frac1{2m-1},$$
which implies that
$$-{2^{-\frac{1}{2m}}}<u_p\le0.$$
Taking into account these bounds on $p$ and $u_p$ and employing notation introduced in the proof of Lemma~\ref{lem:equiv}, one can see that
\begin{multline*}
\partial_\th\De_{1,p,m}(\th,u_p)|_{\th=\pi/4}
=\frac{2^{1-1/(2m)}(2m-1)(p-p_*)(p-p_{**})}{m\,(1-p)},\quad\text{where} \\
p_{**}:=p_{**}(m):=\frac{2m+1 + \sqrt{4(m-1)(m+2)+1}}{4(2m-1)} > p_*.
\end{multline*}
Because of the assumption $p<p_*$, it follows that
$$\partial_\th\De_{1,p,m}(\th,u_p)|_{\th=\pi/4}>0.$$
On the other hand, 
$$\partial_\th\De_{1,p,m}(\th,u)|_{\th=\pi/4}=
\frac{m\,2^{1/m-1}}{(1-p)p}\,
\partial_\th\De_{p,m}(\th,-u-2^{1-1/(2m)}\,p)|_{\th=\pi/4}$$
for all real $u$; this follows from \eqref{eq:De1,De2}, in view of the fact that the derivatives of 
$\cos^{1-1/m}\th\,\sin^{1-1/m}\th$ and $\cos^{1/m}\th+\sin^{1/m}\th$ in $\th$ at $\th=\pi/4$ are zero. 
Hence,
$$\partial_\th\De_{p,m}(\th,t_p)|_{\th=\pi/4}>0
\quad\text{for}\quad t_p:=-u_p-2^{1-1/(2m)}\,p.$$ 
Note also that $\De_{p,m}(\pi/4,t)=0$ for all real $t$. 
Therefore, $\De_{p,m}(\th,t_p)<0$ for all $\th$ in a left neighborhood of 
$\pi/4$. Now \eqref{eq:de} implies that $\pi/4$ is not a point of maximum in 
$\th$ of 
$\E f_t(\BC_1\,\cos^{1/m}\th+\BC_2\,\sin^{1/m}\th)$ for $t=t_p$. 
Hence, in view of  
the homogeneity argument used in the proof of 
Lemma~\ref{lem:equiv},  
$\pi/4$ is not a point of maximum in 
$\th$ of 
$\E f_t(\BS_1\,\cos^{1/m}\th+\BS_2\,\sin^{1/m}\th)$ for 
$t=t_p/\sqrt{pq}$. 
But, for any $t\in\R$, one has $f_t\in\H3\subseteq\F3$ (where the set inclusion follows by Lemma~\ref{lem:F-al}).
Thus, one obtains a contradiction with the assumed statement 
\eqref{it:Bern,n=2} of Theorem~\ref{th:main}. 
\end{proof}



\end{document}